\documentclass[journal]{IEEEtran}

\usepackage{cite}
\usepackage{amsmath,amssymb,amsfonts}
\usepackage{algorithmic}
\usepackage{graphicx}
\usepackage{textcomp}
\usepackage{graphicx,color}
\usepackage{mathrsfs}
\usepackage{subfigure}
\usepackage{url}
\usepackage{booktabs}
\usepackage{array}
\usepackage{mathtools} 		
\usepackage{epstopdf}
\usepackage[vlined,ruled]{algorithm2e}
\usepackage{upgreek}
\usepackage{dsfont}
\usepackage{stackengine} 
\usepackage{graphicx}
\usepackage{subfigure}
\usepackage{caption}
\newcommand{\calX}{{\mathcal{X}}}
\newcommand{\calO}{{\mathcal{O}}}
\newcommand{\calT}{{\mathcal{T}}}

\newcommand{\calN}{{\mathcal{N}}}
\newcommand{\calL}{{\mathcal{L}}}

\newcommand{\bbR}{{\mathbb{R}}}
\newcommand{\bbE}{{\mathbb{E}}}

\newcommand{\bbN}{{\mathbb{N}}}
\newcommand{\bx}{{\mathbf{x}}}
\newcommand{\bz}{{\mathbf{z}}}

\newcommand{\bs}{{\mathbf{s}}}

\newcommand{\be}{{\mathbf{e}}}
\newcommand{\bg}{{\mathbf{g}}}

\newcommand{\by}{{\mathbf{y}}}
\newcommand{\bw}{{\mathbf{w}}}
\newcommand{\bA}{{\mathbf{A}}}
\newcommand{\ba}{{\mathbf{a}}}
\newcommand{\bB}{{\mathbf{B}}}
\newcommand{\bD}{{\mathbf{D}}}

\newcommand{\bK}{{\mathbf{K}}}

\newcommand{\bd}{{\mathbf{d}}}
\newcommand{\bmu}{{\utwi{\mu}}}
\newcommand{\utwi}[1]{\mbox{\boldmath $#1$}}

\newcommand{\blambda}{{\utwi{\lambda}}}
\newcommand{\bkappa}{{\utwi{\kappa}}}
\newcommand{\norm}[1]{\lVert #1 \rVert}

\newtheorem{lemma}{Lemma}
\newtheorem{theorem}{Theorem}

\newtheorem{assumption}{Assumption}
\newtheorem{remark}{Remark}
\usepackage[ruled,vlined]{algorithm2e}

\def\BibTeX{{\rm B\kern-.05em{\sc i\kern-.025em b}\kern-.08em
    T\kern-.1667em\lower.7ex\hbox{E}\kern-.125emX}}

\title{Time-Varying Optimization of Networked Systems with Human Preferences }

\author{Ana M. Ospina, Andrea Simonetto,  Emiliano Dall'Anese
\thanks{A. Ospina and E. Dall'Anese are with the Department of Electrical, Computer and Energy Engineering, University of Colorado Boulder, Boulder, CO, USA; e-mails: \{ana.ospina, emiliano.dallanese\}@colorado.edu. A. Simonetto is with the UMA, {ENSTA Paris}, Institut Polytechnique de Paris, 91120 Palaiseau, France; e-mail: andrea.simonetto@ensta-paris.fr.}
\thanks{The work of A. Ospina and E. Dall'Anese was supported by the National Science Foundation (NSF) CAREER award 1941896, NSF ERC ASPIRE, and by the U.S Department of Energy project 3.2.6.80.}}

\begin{document}
\maketitle

\begin{abstract}
This paper considers a time-varying optimization problem associated with a network of systems, with each of the systems shared by (and affecting) a number of individuals. The objective is to minimize  cost functions associated with the individuals' preferences, which are unknown, subject to time-varying constraints that capture physical or operational limits of the network. To this end, the paper  develops a distributed  online optimization algorithm with concurrent learning of the cost functions. The cost functions are learned on-the-fly based on the users' feedback (provided at irregular intervals) by leveraging tools from shape-constrained Gaussian Process. The online algorithm is based on a primal-dual method, and acts effectively in a closed-loop fashion where: i) users' feedback is utilized to estimate the cost, and ii) measurements from the network are utilized in the algorithmic steps to bypass the need for sensing of (unknown) exogenous inputs of the network. The performance of the algorithm is analyzed in terms of dynamic network regret and constraint violation. Numerical examples are presented in the context of real-time optimization of distributed energy resources.
\end{abstract}


\section{Introduction}


We consider optimization problems associated with a network of systems, where each of the systems is shared by a number of individuals. Typically, such a multiuser problem includes a cost given by a sum of user-specific functions, and a set of constraints that capture physical or operational limits of the network (and, thus, that couple the users' decisions). For example, multiuser problems were considered in \cite{koshal2011multiuser}, and solved via primal-dual  method based on a regularization of the Lagrangian function; and, a distributed resource allocation problem was investigated in \cite{turan2020resilient} (where the communications between users' are susceptible to adversarial attacks). 
Similar formulations arise in a time-varying setting, where the cost and/or constraints evolve over time to reflect changes in the constraints or problem inputs~\cite{popkov2005gradient,SimonettoGlobalsip2014,fazlyab2016self,simonetto2020time}. Specifically, unconstrained problems with a known time-varying cost are analyzed in \cite{popkov2005gradient} via a gradient descent method. A time-varying multiuser problem is presented in \cite{SimonettoGlobalsip2014}, where a double regularization both in the primal and in the dual space is employed to increase the convergence rate of the algorithm. In \cite{fazlyab2016self}, convex optimization problems with known time-varying objective functions are solved using a real-time self-triggered control method. We refer the reader to the survey in \cite{simonetto2020time}, for a complete list of references on time-varying convex optimization.

\begin{figure}[t!] 
    \centering
\includegraphics[scale=0.45]{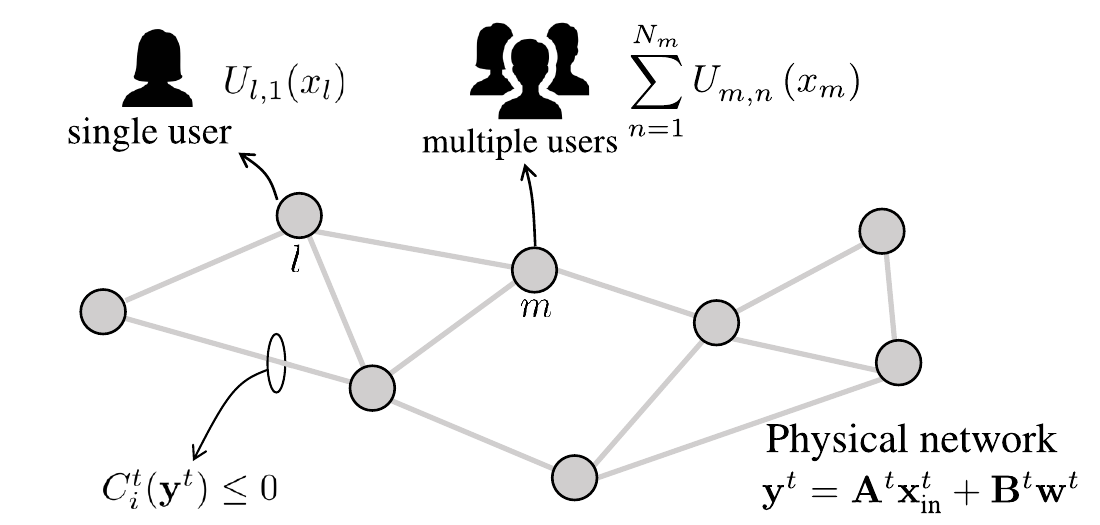}
   \caption{We consider a network of systems, coupled through physical or operational dependencies, where each system is shared by (and affects) a number of users.} 
   \vspace{-.5cm}
    \label{fig:toy_eg}
\end{figure}

Generally, one prerequisite for solving these problems is that the user-specific functions are known or properly crafted based on synthetic or average models. Instead, we consider a case where synthetic costs are not representative of the preferences of individual users, or fail to capture the diversity in their perception of comfort, safety, or dissatisfaction~\cite{bajcsy2021analyzing,LASSEN2020smiley-face,lepri2017tyranny}. More precisely, a learning method that leverage human feedback is explored in \cite{bajcsy2021analyzing}, where a robot observes data in the form of state-action and the discrete dynamics are learned during the execution of the algorithm. An example of field research of occupant satisfaction to the indoor temperature via recording of smiley-face polling station is presented in \cite{LASSEN2020smiley-face}. 
Moreover, \cite{lepri2017tyranny} presents an overview of opportunities and risks of data-driven decision-making methods.

We tackle the problem of solving a network optimization problem when the users' cost functions are  \emph{unknown}, and develop an  online algorithm where the cost functions are learned concurrently with the execution of the algorithmic steps. In particular, the learning procedure leverages ``users' feedback'', and utilizes learning tools from shape-constrained Gaussian Processes (GPs) \cite{jrnl_wang}. The algorithm is implemented in a \emph{distributed} fashion; this allows users sharing a system to agree on a solution that minimizes the sum of their (learned) functions, without revealing their preferences or their feedback. In this work, we consider the case where the users report their preferences trustworthy (and we do not consider adversarial behaviors). In addition to the feedback from the users, the algorithm leverages the ideas from, e.g.,~\cite{bolognani2013distributed, Dallanese2019_PD}, and utilizes measurements of  network outputs instead of the network model.

The endeavors are motivated by a number of problems arising in power systems \cite{wang2014adaptive,Lesage_2020,bolognani2013distributed}, charging of electric vehicles \cite{EV_eg1}, and human-aware robot systems \cite{luo2020sociallyaware}, just to name a few. While Section~\ref{sec:problemstatement} will explain some examples, we stress here that an accurate model of the comfort/satisfaction level of users is typically unknown, and it may vary not only across individuals, but also for the same individual~\cite{bajcsy2021analyzing,LASSEN2020smiley-face}. 
In this context, shape-constrained GP is a powerful tool for  nonparametric  function  estimation  when  the  underlying function is presumed to be convex or monotonic~\cite{jrnl_wang}; this  implies  that  each  user  has  a  preference for a  set of solutions. 

\emph{Contributions}. Overall, the main contributions of the paper are as follows. C1) We show how to estimate the users' functions via GPs and in particular,  how to handle  asynchronous  and noisy data easily. C2) We develop an  online  algorithm  based  on  a primal-dual  method to solve the formulated time-varying network optimization problem; the algorithm is implemented online, and it is modified to accommodate measurements and learned functions. C3) We propose a consensus-based formulation that leads to a distributed online algorithm where users minimize the sum of their (learned) functions, without revealing their preferences or their feedback. C4) For the performance analysis, we view the algorithm  as an inexact online primal-dual method where errors are due to the shape-constrained GP approximation and errors in the computation of the gradient of the GP. We derive bounds for a dynamic extension of the network regret~\cite{Online_PD,lee2016distributed}, and for the average constraint violation (see also, e.g.,~\cite{chen2018bandit,hall2015online,Shahrampour2018,yi2016tracking} and references therein for other notions of dynamic regret). And, C5) we showcase the proposed methodology in an example related to control of distributed energy resources in power grids~\cite{Lesage_2020}. 

\emph{Prior works}. In the context of bandit optimization, a static regret analysis is performed in \cite{chen2018bandit}, where the objective is to minimize a sequence of unknown convex cost functions and one has access to functional evaluations. Time-varying cost and constraints are considered in the context of a bandit setting in \cite{flaxman2004online} and online saddle-point algorithms are utilized, where the gradient information is acquired via one-point or multi-point estimates. In \cite{li2020online}, a gradient-free approach is presented to solve a multi-agent distributed constrained optimization problem where the goal is to cooperatively minimize the sum of time-changing local cost functions subject to time-varying coupled constraints. As an example of works in the context of  zeroth-order methods, a multi-agent optimization problem is analyzed in \cite{tang2020zeroth} where the objective is to minimize the average of the nonconvex local costs that depends on the joint actions of the agents. See also \cite{liu2020primer} for an overview of zeroth-order algorithms. Our approach is different from zeroth-order or bandit methods~\cite{flaxman2004online,chen2018bandit,liu2020primer,tang2020zeroth,li2020online}; these would require multiple functional evaluations at each step to estimate the gradient of the users' costs -- something not feasible for our problem. Our contribution is to use shape-constrained GPs to process feedback that is provided parsimoniously and at irregular intervals. The considered setting is also different from~\cite{Srinivas2012}, where GPs are utilized to maximize an unknown function. 

A similar problem setup was studied in~\cite{notarnicola2020distri}; however, no constraints are considered (in fact, a gradient tracking scheme is utilized) and the user's cost is assumed to have a known quadratic structure a priori. We significantly extend our previous works~\cite{GP_userfeedback,OspinaSHGP2020} by considering a distributed setting and constrained problems. 
Relative to the works \cite{Online_PD,lee2016distributed} on distributed online primal-dual algorithms, we consider problems where the optimal solution set changes in time (and it is not fixed in time), and with both consensus constraints and time-varying constraints  that  capture  operational  limits  of  the  network. We also analyze the network regret when the online algorithm is \emph{inexact}. Indeed, we can recover an asymptotic bound of $\calO(T^{\frac{1}{2}})$  as in~\cite{Online_PD,lee2016distributed} from our  regret bound when the cost functions are known, the gradients are computed exactly, and the optimal solutions are not time-varying (with $T$ the number of time steps). Our analysis can also be utilized to extend~\cite{tang2020zeroth} to time-varying settings. Finally, we extend the findings of \cite{Dallanese2019_PD,bolognani2013distributed} by considering consensus constraints (and accompanying distributed implementations), and by including shape-constrained GPs in the overall algorithmic framework.

The remainder of this paper is organized as follows. Section II presents the time-varying problem, and Section III outlines the GP-based learning approach and the proposed online algorithm. Section IV presents the performance analysis. The numerical results for a demand response problem are in Section V, and Section VI concludes the paper\footnote{\textit{Notation:} Upper-case (lower-case) boldface letters will be used for matrices (column vectors); $(\cdot)^\top$ denotes transposition. For a given column vector $\bx \in \bbR^n$, $\norm{\bx} := \sqrt{\bx^\top\bx}$; we denote with $\mathbf{1}$ and $\mathbf{0}$ the vectors of all ones and zeros, respectively, where the dimension will be clear from the context. $\mathcal{O}(\cdot)$ refers to the big-O notation, whereas $o(\cdot)$ refers to the little-o notation. The max operator is defined as $[x]^+ := \max\{0, x\}$. The diameter of $\calX^t$ is defined as $\mathrm{diam}(\calX^t) = \max\{ \|\bx-\bx'\| \, : \, \bx,\bx' \in \calX^t \} $. For a given random variable $\xi \in \mathbb{R}$, $\bbE[\xi]$ denotes the expected value of $\xi$.}.

\section{Problem Statement}
\label{sec:problemstatement}

We consider a network of $M$ systems or devices, with each system interacting with (or shared by) a number of users as illustrated in Figure~\ref{fig:toy_eg}. In particular, we denote  by $\calN_m = \{1, \dots, N_m \}$ the set of $N_m \geq 1$ users interacting with the $m$th system. The steady-state input-output relationship of the network is modeled through the time-varying map
\begin{equation}
\label{eq:system}
    \by^t = \bA^t \bx^t_{\textrm{in}} + \bB^t \bw^t
\end{equation}
where $t$ is the time index; $\bx_{\textrm{in}}^t := [x_1^t, \ldots, x_M^t]^\top$, with $x^t_m \in \mathbb{R}$ the controllable input of the $m$th system; $\by^t \in \mathbb{R}^Y$ is a vector of observables or outputs of the network; $\bw^{t} \in \mathbb{R}^W$ is a vector of \emph{unknown} exogenous inputs; and, $\mathbf{A}^t$ and $\mathbf{B}^t $ are possibly time-varying matrices (with $\bA^t$ known and $\bA^t \neq \mathbf{0}$) of appropriate dimensions. We assume that $t \in \calT := \{k \Delta, k \in \bbN \}$, where $\Delta > 0$ is a given time interval. 

Let $U_{m,n}: \mathbb{R} \rightarrow \mathbb{R}$ be a convex function representing the cost of a user $n \in \calN_m$ associate with the $m$th device\footnote{ We note that the problem formulation and solution approach proposed in this paper can be naturally extended to the case where the number of users interacting with a device changes over time.}; this function can capture a sense of dissatisfaction, discomfort, or simply preferences  depending on the specific application (a few examples will be provided shortly). The goal is then to compute the sequence of inputs that solves the following time-varying optimization problem~\cite{popkov2005gradient, SimonettoGlobalsip2014,fazlyab2016self} associated with the network:  
\begin{subequations}
\label{eq:P0}
\begin{align}
\underset{\{x_m \in \calX_m^t \}_{m=1}^M}{\min} &\quad \sum_{m=1}^M \sum_{n=1}^{N_m} U_{m,n} \left(x_{m}\right) \\
  \text{s.~to: } & C^t_i \left(\bA^t \bx_{\textrm{in}} + \bB^t \bw^t \right) \leq 0, i = 1, \dots, N_C  \label{eq:P0_c1}
\end{align}
\end{subequations}
where $\{\calX_m^t\}_{m=1}^M$ are time-varying convex and compact sets for the inputs, and $C^t_i: \bbR^M \rightarrow \bbR$ is a time-varying  function parameterized by $\bw^t$. In particular,~\eqref{eq:P0_c1} models operational constraints associated with the network output $\by^t$. The main goal of~\eqref{eq:P0} is to generate a sequence of inputs that minimizes the cost of the users, while respecting pertinent network constraints. To this end, this paper focuses on addressing the following four main challenges related to~\eqref{eq:P0}: 

\noindent $\bullet$ \emph{Challenge~1:} The variable $x_m$ is \emph{coupled} across the users $\calN_m$. Our goal is to solve~\eqref{eq:P0} in a distributed fashion, where users are not required to share information about their costs $\{U_{m,n}\}$ (and their functional evaluations) with the device or the network operator.

\noindent $\bullet$ \emph{Challenge~2:} The costs $\{U_{m,n}\}$ are \emph{unknown}. One may utilize synthetic cost functions for users’ preferences, comfort or satisfaction based on statistical models; however, these synthetic costs may fail to capture the  diversity in the users' perception, preferences, and goals. 

\noindent $\bullet$ \emph{Challenge~3:} The exogenous inputs $\bw^{t}$ may be  \emph{unknown} or partially known. 

\noindent $\bullet$ \emph{Challenge~4:} The constraints $C_i^t$ and the sets $\{\calX_m^t\}$ may change at each time $t$; the interval $\Delta$ may not be sufficient to solve each instance of the time-varying problem~\eqref{eq:P0} to convergence.  

To address the challenges 1-4 above, this paper will develop a \emph{consensus-based} \emph{online} algorithm where the cost functions of the users are \emph{learned} concurrently with the execution of the algorithm, and \emph{measurements} of $\by^t$ are utilized in the algorithmic updates.

Before proceeding, we list representative examples of applications in networked systems in which the problem~\eqref{eq:P0} and the accompanying online optimization framework considered in this paper are particularly well-suited. 

\textbf{\emph{Example~1}} \emph{(Power grids)}. In the context of power grids,~\eqref{eq:P0} may capture a demand response task~\cite{wang2014adaptive,Lesage_2020} or a real-time optimal power flow problem~\cite{bolognani2013distributed}. In this case, $\bx_\textrm{in}$ are the power setpoints of distributed energy resources, $\bw^t$ is a vector of powers consumed by uncontrollable loads, and $y^t$ may represent the total power as $y^t = \mathbf{1}^\top \bx_{\textrm{in}} + \mathbf{1}^\top \bw^t$ and/or voltage magnitudes (one can also compute the matrix $\mathbf{A}^t$ based on a linearized AC model). Constraints~\eqref{eq:P0_c1} may ensure that the net power follows a given automatic gain control or demand response signal; e.g., $C_i^t(y) = (y - y_{\textrm{ref}}^t)^2 - \zeta$, where $y_{\textrm{ref}}^t$ is a given reference signal and $\zeta > 0$ is an error tolerance. One may also have voltage constraints. The costs $\{U_{m,n}\}$ model discomfort, e.g., the indoor temperature. The number of users depends on the particular setting, for example one may expect a large number of users sharing a conference room in an office building. The users that share a device $m$ agrees on the input $x_m$; for example, people sharing a room agree on a temperature and, thus, on the use of an HVAC system.   \hfill $\Box$


\textbf{\emph{Example~2}} \emph{(Electric vehicle (EV) charging)}. The formulation~\eqref{eq:P0} may represent the problem of charging a fleet of vehicles at a charging station; the variable $x_m$ represents the charging rate (or the expected time for full charging) of a vehicle, constraint~\eqref{eq:P0_c1} captures limits on the total power consumed by the charging station, while $\{ U_{m,n}\}$ capture  the dissatisfaction for a given expected time of charging completion~\cite{EV_eg1}.  If a EV is shared by multiple users, $x_m$ represents the agreed expected charging time $x_m$. \hfill $\Box$

Additional examples in transportation systems, communication systems, and human-aware robot systems are not included for space limitations.   

We now proceed with the reformulation of~\eqref{eq:P0} into a consensus-based problem. To this end, let $x_{m,n}$ be an auxiliary variable representing the input for the $m$th system preferred by user $n \in \calN_m$; with the auxiliary variables in place, stack the optimization variables in the column vector
$$\bx^t = \left[x_{1}^t,  x_{1,1}^t, \dots, x_{1, N_1}^t, \dots,  x_{M}^t, x_{M,1}^t, \dots,  x_{M, N_M}^t \right]^\top.$$ 
Accordingly,~\eqref{eq:system} can be rewritten as $ \by^t = \bar{\bA}^t \bx^t + \bB^t \bw^t$, where $\bar{\bA}^t$ is an augmented version of $\bA^t$ with appropriate zero entries.

Then,~\eqref{eq:P0} can be equivalently reformulated as:  
\begin{subequations}
\label{eq:P1}
\begin{align}
\underset{\{x_m, \, x_{m,n} \in \calX_m^t \}} {\min} &\quad  \sum_{m=1}^M \sum_{n=1}^{N_m} U_{m,n} \left(x_{m,n}\right) \\ 
 & \hspace{-1.4cm} \text{s. to: }  C^t_i \left(\bar{\bA}^t \bx + \bB^t \bw^t\right) \leq 0 , i = 1, \dots, N_C  \hspace{-.2cm} \label{eq:P1_c1} \\ 
  & \hspace{-.5cm} x_{m} = x_{m,n} , \,  \forall \, n \in \calN_m, \, \forall \,  m = 1, \dots, M \label{eq:P1_c2}
\end{align}
\end{subequations}
where~\eqref{eq:P1_c2} defines consensus constraints for the $m$th system and its users. We note that the choice of posing the consensus constraints as in~\eqref{eq:P1_c2} is not generic; the particular structure of~\eqref{eq:P1_c2}, which can be modeled as a star graph (with the system as the central node) will facilitate the development of a closed-loop algorithm where $\{x_{m}\}$ will be physical inputs (as in~\eqref{eq:system}) and $\{x_{m,n}\}$ will be auxiliary control variables.         

For future developments, define $f(\bx) := \sum_{m=1}^M \sum_{n=1}^{N_m} U_{m,n} \left(x_{m,n}\right)$ and  $\mathcal{X}^t:= \mathcal{X}^t_1 \times \dots \times \mathcal{X}^t_M \subseteq \bbR^M$. For each device $m$, define  the  incidence matrix $\bD_m \in \bbR^{N_m \times (1+N_m)}$, which represents the relation between the device and the users sharing that device. For example, if the device $1$ is shared by 2 users we have that
\begin{equation*} \bD_1 =
    \begin{bmatrix}
    1 & -1 & 0 \\
    1 & 0 & -1
    \end{bmatrix}.
\end{equation*}
By organizing these incidence matrices in a block diagonal matrix, we construct the augmented incidence matrix $\bD = \mathrm{diag}(\bD_1, \bD_2, \dots, \bD_M)$. Then, \eqref{eq:P1_c2} can be compactly written as $\bD \bx = \mathbf{0}$, capturing  consensus constraints over $M$ systems.

The next section outlines the proposed algorithm to solve the time-varying network optimization problem~\eqref{eq:P0}. We note that in what follows we consider the case where $N_C = 1$ to simplify exposition and notation; however, the technical arguments of this paper straightforwardly  extend to problems with multiple constraints on $\by^t$.

\section{Distributed Online Algorithm} 
\label{sec:distri_PD}

In this section, we present our \emph{consensus-based} \emph{online} algorithm for solving~\eqref{eq:P1}. We first explain how to derive a distributed online algorithm based on a primal-dual method; then, we elaborate on how to estimate the cost functions based on feedback received from the users at infrequent times.

\subsection{Online primal-dual algorithm}

We start by defining the time-varying Lagrangian function associated with~\eqref{eq:P1} as follows: 
\begin{align}
    \calL^t(\bx^t, \nu^t, \blambda^t) : = & f(\bx^t) + \nu^t \, C^t \left( \bar{\bA}^t \bx^t + \bB^t \bw^t \right) + (\bD \bx^t)^\top \blambda^{t} \nonumber
\end{align}
where $\nu^t \in \bbR_+$, $\blambda^t \in \bbR^N$ are the dual variables associated with the constraint~\eqref{eq:P1_c1}, and the $M$ consensus constraints~\eqref{eq:P1_c2}, respectively. 

Recall that $t \in \calT := \{k \Delta, k \in \bbN \}$, where $\Delta > 0$ represents the  time required to perform one step of the online algorithm (and with the constraints and the vector $\bw^t$ possibly changing at every interval $\Delta$). For a given step size $\alpha > 0$, a \emph{model-based} online projected primal-dual algorithm involves the sequential execution of the following steps:
\begin{subequations}
\label{eq:updates}
    \begin{align} 
        \bx^{t+1} &= \mathrm{proj}_{{\calX^t}} \left\{\bx^t -\alpha \left( \nabla_{\bx} f(\bx^t) \right. \right. \nonumber \\
        & \left. \left. \,\,\,\, + \nu^t (\bar{\bA}^t)^\top \nabla_{\bx} C^t(\bar{\bA}^t \bx^t + \bB^t \bw^t) + \bD^\top \blambda^t \right) \right\}, \label{eq:update_x} \\ 
        \nu^{t+1} &= \mathrm{proj}_{\Psi^t} \left\{\nu^t + \alpha C^t \left( \bar{\bA}^t  \bx^t + \bB^t \bw^t \right) \right\}, \label{eq:update_nu} \\ 
        \blambda^{t+1} &= \mathrm{proj}_{{\Lambda^t}} \left\{ \blambda^t + \alpha \bD \bx^t \right\} , \label{eq:update_la}
    \end{align}
\end{subequations}
where $\Psi^t$ is a convex and compact set, and $\Lambda^t = \Lambda_{1,1}^t \times \dots \times \Lambda_{1,N_1}^t \times \dots \times \Lambda_{M,N_M}^t$, with $ \Lambda_{m,n}^t$ also a convex and compact set constructed as explained shortly in Section~\ref{ref:performance} (and in, e.g.,~\cite{koshal2011multiuser} and~\cite{Online_PD}). 

We note that the steps~\eqref{eq:update_x} and~\eqref{eq:update_nu} require one to know $\bw^t$ and $\bB^t$ (as explained in Challenges 3 and 4). To bypass this hurdle, we adopt the strategy proposed in, e.g.,~\cite{bolognani2013distributed,Dallanese2019_PD}. To this end, assume that $\bx_{\mathrm{in}}^t$ are implemented as inputs to the systems at time $t$, and let $\hat{\by}^t$ represent a \emph{measurement} of $\by^t$ collected at time $t$. Then, $\nu^{t+1}$ can be computed as $\nu^{t+1} = \mathrm{proj}_{\Psi} \left\{\nu^t + \alpha C^t \left( \hat{\by}^t \right) \right\}$; it is clear that information about $\bw^t$ and $\bB^t$ is not required in this case. We also note that the step~\eqref{eq:update_x} would require the gradient of the functions $\{U_{m,n}\}$, which are not known (see Challenge 1). Let $\hat{U}^t_{m,n}(x^t_{m})$ represent an estimate of the function $U_{m,n}(x^t_{m})$  available at $x^t_{m}$, and let $g_{m,n}^t$ denote an \emph{estimate} of the gradient (or derivative, in this case) of $U_{m,n}(x^t_{m})$ at time $t$. With these definitions, the steps in~\eqref{eq:updates} can be suitably modified to accommodate measurements of $\by^t$ and estimates of the gradient of the cost function: 
\begin{subequations}
\label{eq:OPD_case1}
\begin{align} 
    x_m^{t+1} &= \mathrm{proj}_{{\calX_m^t}} \left\{ x_m^t -\alpha \left( \nu^t (\ba_m^t)^\top  \nabla C^t(\hat{\by}^t) + \sum_{n=1}^{N_m} \lambda^t_{m,n} \right) \right\}  \label{eq:OPD_sys} \\ 
    x_{m,n}^{t+1} &= \mathrm{proj}_{{\calX_m^t}} \left\{  x_{m,n}^t -\alpha \left( g_{m,n}^t - \lambda^t_{m,n} \right) \right\} \;  \forall \, n \in \calN_m, \, \forall \,  m,  \label{eq:OPD_user} \\ 
    \nu^{t+1} &= \mathrm{proj}_{\Psi^t} \left\{\nu^t + \alpha C^t \left( \hat{\by}^t \right) \right\}, \label{eq:OPD_nu} \\ 
    \lambda_{m,n}^{t+1} & = \mathrm{proj}_{{\Lambda_{m,n}^t}} \left\{ \lambda_{m,n}^t + \alpha \left(x_m^t - x_{m,n}^t \right) \right\} ,  \forall \, n \in \calN_m, \, \forall \,  m,  \label{eq:OPD_lambda}
\end{align}
\end{subequations}
where $\ba_m^t$ is the $m$th column of the matrix $\bA^t$, and we note that step~\eqref{eq:OPD_sys} is performed in parallel at each system.  While measurements of $\by^t$ are utilized in the primal step~\eqref{eq:OPD_sys} and in the dual step~\eqref{eq:OPD_nu}, an estimate of the gradient $g_{m,n}^t$ will be obtained by leveraging feedback from the users as explained in the following subsection. Subsequently, Section~\ref{sec:algorithm} will overview the steps of the proposed algorithm.   

\subsection{Online Learning via Shape-Constrained GPs}
\label{eq:shape-constrained}

We assume that the feedback of the user comes in the form of a (possibly noisy) functional evaluation; that is, for a given point $x$, one may receive from a user $n$ interacting with the device $m$ a rating $z$ given by $z = U_{m,n}(x) + \epsilon$. With this model, a popular way to obtain the gradient of $U_{m,n}$ at a point $x$ is via zeroth-order methods (see, e.g.,~\cite{flaxman2004online,liu2020primer} and references therein). However, zeroth-order methods are not well-suited for our algorithm, because they would require feedback from each individual user at \emph{each iteration} of the algorithm (either a single functional evaluation or more, depending on the particular method utilized). Instead, we consider a more realistic case where users may provide feedback more parsimoniously and at irregular intervals. 

With this setting in mind, we consider utilizing the noisy functional evaluations provided by a user to estimate the function $U_{m,n}$ via GPs~\cite{GPML_Rasmussen06}. In particular, we apply a GP with specific constraints on the shape of the function. As explained in~\cite{jrnl_wang}, shape-constrained GPs offers a powerful tool for nonparametric function estimation  when the underlying function is presumed to be convex or monotonic. In our case, convexity  or monotonicity implies that each user has a preference for a particular solution (or a convex set of solutions), and the dissatisfaction or discomfort increases as the system deviates from such a preferred point.  In the following, we first briefly introduce GPs, and then we provide the main equations to estimate $U_{m,n}$  via a shape-constrained GPs. Since we utilize a GP per function $U_{m,n}$, we drop the  subscripts $m$, $n$ for simplicity of exposition. We also note that the extension to the multi-dimensional shape-constrained GP is possible (see, for example,~\cite{jrnl_wang}); however, it is left as future work.

A GP is a stochastic process $U(x)$ and it is specified by its mean function $\mu({x})$ and its covariance function $k(x, x')$; where $k(\cdot, \cdot)$ is a given kernel; i.e., for any $x, x' \in \mathcal{X} \subseteq \mathbb{R}$, $\mu({x}) = \mathbb{E}[U({x})]$ and $k({x}, {x}') = \mathbb{E}[(U({x}) - \mu({x}))(U({x}') - \mu({x}'))]$ \cite{GPML_Rasmussen06}. Hereafter, we use the short-hand notation $U(x) \sim \mathcal{GP} (\mu(x), k(x,x'))$ for $U({x})$ defined as a GP with mean $\mu(x)$ and covariance $k(x,x')$. Let $\mathbf{x}_p = [ x_1 \in \mathcal{X}, \dots, x_p \in \mathcal{X}]^\top$ be the set of $p$ sample points collected over a given period\footnote{In the context of our work, the vector $\mathbf{x}_p$ associated with a system $m$ and a user $n$ would contain the points $x_{m,n}^{t_1}, \ldots, x_{m,n}^{t_p}$, where $\{t_i\}_{i = 1}^p \subset \mathcal{T}$ are the times instants where the user provides feedback.}; let ${z}_i = U({x}_i) + \epsilon_i$, with $\epsilon_i \stackrel{\mathrm{iid}}{\sim} \mathcal{N}({0}, \sigma^2)$ Gaussian noise, be the noisy measurement at the sample point $x_i$; define further $\mathbf{z} = [z_1, \dots, z_p]^\top$. Then, the posterior distribution of $(U(x)|\bx_p, \bz)$ is a GP  with mean $\mu({x})$, covariance $k({x}, {x}')$, and variance $\sigma^2(x)$ given by~\cite{GPML_Rasmussen06}:
\begin{subequations} \label{eq:GP_update}
\begin{align}
    \mu({x}) &= \mathbf{k}(x)^\top (\mathbf{K} + \sigma^2 \mathbf{I})^{-1} \mathbf{z}, \\
    k({x}, {x}') &= k(x,x') - \mathbf{k}(x)^\top (\mathbf{K} + \sigma^2 \mathbf{I})^{-1} \mathbf{k}(x'), \\
    \sigma^2(x) &= k(x,x),
\end{align}
\end{subequations}
where $\mathbf{k}(x) = [k(x_1,x), \dots, k(x_p,x)]^\top$, $\mathbf{K}$ is the positive definite kernel matrix $[k(x,x')]$. Note that \eqref{eq:GP_update} depends on the $p$ sample points for a specific period of time. Thus, an estimate of the unknown function $U({x})$ can be written as $U({x}) \sim \mathcal{GP}( \mu({x}),  k({x}, {x}'))$. Using a squared exponential (SE) kernel as an example, $k({x},{x}')$ is given by $k({x},{x}') = \sigma_{f}^2 \mathrm{e}^{-\frac{1}{2 \ell ^2}({x}-{x}')^{2}}$ for the univariante input case, where the hyperparameters are the signal variance $\sigma_{f}^2$ and the characteristic length-scale $\ell$.

We now explain how to enforce shape constraints through the use of derivative processes (noting that since differentiation is a linear operator, derivatives of a GP remain a GP \cite{GPML_Rasmussen06}). Because the GP obtained from the SE kernel function has derivatives of all orders, the mean and covariance function (jointly with the original process and the second-order derivative process) are \cite{jrnl_wang}:

\begin{subequations}
\label{eq:GP_dsecond}
\begin{equation}
   \hspace{-4.0cm} \mathbb{E}\left[\frac{\partial^2 U({x})}{\partial {x}^2}\right] = \frac{\partial^2 \mu({x})}{\partial {x}^2} = {0},
\end{equation}
\vspace{-0.5cm}
\begin{multline}
        k^{22}(x,x') := \text{cov}\left[\frac{\partial^2 U({x})}{\partial {x}^2}, \frac{\partial^2 U({x}')}{\partial {x}'^2}\right] = \sigma_{f}^2  \mathrm{e}^{-\frac{1}{2\ell^2}(x-x')^{2} } \times \\ \frac{1}{\ell^4}\left(\frac{1}{\ell^4}(x-x')^4 - \frac{1}{\ell^2}6(x-x')^2+3\right), 
\end{multline}
\vspace{-0.5cm}
\begin{multline}
k^{02}(x,x') :=    \text{cov}\left[\frac{\partial^2 U({x})}{\partial {x}^2}, U({x}') \right]  \\ = \sigma_{f}^2 \mathrm{e}^{ -\frac{1}{2\ell^2}(x-x')^{2} }  \left(\frac{1}{\ell^4}(x-x')^2 - \frac{1}{\ell^2}\right).
\end{multline}
\end{subequations}
We consider enforcing constraints on the second-order derivative at $q$ points $\mathbf{s} := [s_1, \dots, s_q]^\top$ (and these points may be in general different from $\bx_p$). Let $\mathbf{U}(\mathbf{x}_p)= [U(x_1), \dots, U(x_p)]^\top$ and $\mathbf{U}''(\mathbf{s})= [U''(s_1), \dots, U''(s_q)]^\top$; then, the joint distribution of the GP and its second-order derivative is:
\begin{equation*}
    \begin{bmatrix} 
    \mathbf{U}(\mathbf{x}_p) \\
    \mathbf{U}''(\mathbf{s})
    \end{bmatrix}
    \sim \mathcal{N} \bigg( \begin{bmatrix} 
    \mu \mathbf{1}_p\\
    \mathbf{0}_q
    \end{bmatrix},
    \begin{bmatrix} 
    \mathbf{K}(\mathbf{x}_p,\mathbf{x}_p) & \mathbf{K}^{\mathbf{02}}(\mathbf{x}_p,\mathbf{s}) \\
    \mathbf{K}^{\mathbf{20}}(\mathbf{s},\mathbf{x}_p) & \mathbf{K}^{\mathbf{22}}(\mathbf{s},\mathbf{s})
    \end{bmatrix} \bigg),
\end{equation*}
where, 
\begin{align*}
    \mathbf{K}(\mathbf{x}_p,\mathbf{x}_p) &= \mathbf{K}, \qquad \qquad \qquad
    \mathbf{K}^{\mathbf{02}}(\mathbf{x}_p,\mathbf{s}) = [k^{02}(x, s)],\\
    \mathbf{K}^{\mathbf{20}}(\mathbf{s},\mathbf{x}_p) &= \mathbf{K}^{\mathbf{02}}(\mathbf{x}_p,\mathbf{s})^\top, \qquad \;
    \mathbf{K}^{\mathbf{22}}(\mathbf{s},\mathbf{s}) = [k^{22}(s, s')].
\end{align*}
Next, assign to $U(\cdot)$ a GP prior, and consider aiming at an estimated function that is $L_U$-smooth and $\gamma_U$-strongly convex, for a given $L_U > 0$ and $\gamma_U > 0$. Then, following \cite[Lemma 3.1]{jrnl_wang}, the joint conditional posterior distribution of $(\mathbf{U}({x}^\circ)|\mathbf{U}''(\mathbf{s}), \mathbf{x}_p, \mathbf{z})$, for a point ${x}^\circ \in \mathbb{R}$, given the current observations $\mathbf{z}$, is a GP with mean, covariance, and standard deviation given by:
\begin{subequations}
\label{eq:shape_GP}
    \begin{align}
        \label{eq:GP_def_mu}
        & \bar{\mu}({x^\circ}) = \mu + B_3(\bx_p, {x}^\circ,\mathbf{s}) B_1(\bx_p,\mathbf{s})^{-1}(\mathbf{z} - \mu \mathbf{1}_p) \nonumber \\
        & + (A_2({x}^\circ,\mathbf{s}) - B_3(\bx_p, {x}^\circ,\mathbf{s}) B_1(\bx_p,\mathbf{s})^{-1} A_1(\bx_p,\mathbf{s}))\mathbf{U}''(\mathbf{s}), \\
        & \bar{k}({x^\circ}, {x^\circ}') = A(\bx_p, {x}^\circ,\mathbf{s}), \\
        \label{eq:GP_def_sigma}
        & \bar{\sigma}(x^\circ) = \sqrt{A(\bx_p, {x}^\circ,\mathbf{s})},
    \end{align}
\end{subequations}
and the posterior distribution of $(\mathbf{U}''(\mathbf{s})| \mathbf{x}_p, \mathbf{z})$ is given by:
\begin{equation*}
    (\mathbf{U}''(\mathbf{s})|\mathbf{x}_p,\mathbf{z}) \propto \mathcal{N}(\bmu(\mathbf{s}), \mathbf{S}(\mathbf{s}, \mathbf{s})) \mathbf{1}_{\{ \gamma_U \leq U''(s_i) \leq  L_U, \, i=1, \dots, q \}}
\end{equation*}
where $\mathbf{1}_{\{\cdot\}}$ is a set indicator function,  $(\mathbf{U}''(\mathbf{s})|\mathbf{x}_p,\mathbf{z})$ is a truncated normal distribution and, 
\begin{subequations}
    \begin{align*}
        \bmu(\mathbf{s}) &= \mathbf{K}^{\mathbf{20}}(\mathbf{s},\mathbf{x}_p)(\sigma^2 \mathbf{I} + \mathbf{K}(\mathbf{x}_p,\mathbf{x}_p))^{-1}(\mathbf{z} - \mu \mathbf{1}_p),\\
        \mathbf{S}(\mathbf{s}, \mathbf{s}) &= \mathbf{K}^{\mathbf{22}}(\mathbf{s},\mathbf{s}) \\
        & \hspace{.3cm} -\mathbf{K}^{\mathbf{20}}(\mathbf{s},\mathbf{x}_p)(\sigma^2 \mathbf{I} + \mathbf{K}(\mathbf{x}_p,\mathbf{x}_p))^{-1}\mathbf{K}^{\mathbf{02}}(\mathbf{x}_p,\mathbf{s}),\\
        A_1(\bx_p, \bd) &= \mathbf{K}^{\mathbf{02}}(\mathbf{x}_p,\mathbf{s})\mathbf{K}^{\mathbf{22}}(\mathbf{s},\mathbf{s})^{-1}, \\
        A_2({x}^\circ,\mathbf{s}) &= \mathbf{K}^{\mathbf{02}}({x}^\circ,\mathbf{s})\mathbf{K}^{\mathbf{22}}(\mathbf{s},\mathbf{s})^{-1},\\
        B_1(\mathbf{x}_p,\mathbf{s}) &= \sigma^2 \mathbf{I} + \mathbf{K}(\mathbf{x}_p,\mathbf{x}_p) \\
        & \hspace{.3cm} - \mathbf{K}^{\mathbf{02}}(\mathbf{x}_p,\mathbf{s})\mathbf{K}^{\mathbf{22}}(\mathbf{s},\mathbf{s})^{-1}\mathbf{K}^{\mathbf{20}}(\mathbf{s},\mathbf{x}_p),\\
        B_2({x}^\circ,\mathbf{s}) &= \mathbf{K}({x}^\circ,{x}^\circ) \\
        &\hspace{.3cm}  - \mathbf{K}^{\mathbf{02}}({x}^\circ,\mathbf{s})\mathbf{K}^{\mathbf{22}}(\mathbf{s},\mathbf{s})^{-1}\mathbf{K}^{\mathbf{20}}(\mathbf{s},{x}^\circ),\\
        B_3(\bx_p, {x}^\circ,\mathbf{s}) &= \mathbf{K}({x}^\circ,\mathbf{x}_p) \\
        & \hspace{.3cm} - \mathbf{K}^{\mathbf{02}}({x}^\circ,\mathbf{s})\mathbf{K}^{\mathbf{22}}(\mathbf{s},\mathbf{s})^{-1}\mathbf{K}^{\mathbf{20}}(\mathbf{s},\mathbf{x}_p) ,\\
        A(\bx_p,{x}^\circ,\mathbf{s})) &= B_2({x}^\circ,\mathbf{s}) \\
        & \hspace{.3cm} -  B_3(\bx_p,{x}^\circ,\mathbf{s}) B_1(\mathbf{x}_p,\mathbf{s})^{-1} B_3(\bx_p,{x}^\circ,\mathbf{s})^\top,
    \end{align*}
\end{subequations}
with $\mu$ and $\sigma^2$ given parameters of the prior. The hyperparameters $\ell$ and $\sigma_f^2$ of the GP can be estimated, for example, by using the maximum likelihood estimator \cite{GPML_Rasmussen06}. On the other hand, $L_U$ and $\gamma_U$ can be estimated via cross-validation. The locations of the virtual derivative points are defined beforehand. Note that the shape-constrained GP is a GP whose second-order derivative is constrained on a set of points. When we consider a GP prior with a  differentiable kernel function, we have the advantage that their derivative processes are also GPs and are jointly Gaussian with the original processes \cite{jrnl_wang}.

In this context, shape-constrained GPs guarantee that the posterior mean function $\bar{\mu}(x^{\circ})$ is practically\footnote{With a limited number of enforcing points $\mathbf{s}$, the mean has the functional properties we require up to a very (almost negligible) small error, which is smaller and smaller increasing the number of data points $\bx_p$. The error will be incorporated into the gradient estimation inaccuracies.} smooth and strongly convex~\cite{jrnl_wang}. We will now incorporate the estimates of the users' function in the algorithm as explained next.

\begin{remark} 
One may want to select a dense set of points $\bs$ that ensures a uniform covering of the domain. However, it is important to note that the number of points $q$ influences the dimensions of the matrices $\bK^{\mathbf{20}}$, $\bK^{\mathbf{02}}$, and $\bK^{\mathbf{22}}$; therefore, $q$ should be selected to ensure a sufficiently dense covering, but based on computational complexity considerations. 
\end{remark}

\begin{remark}
In this paper, we use the SE kernel. In this case, two hyperparameters are required for $U_{m,n}$: $\ell$ and $\sigma_f^2$. These two hyperparameters are learned by maximizing the likelihood \cite{GPML_Rasmussen06}. 
Note that the ability to impose shape constraints on the posterior GP is not affected by the kernel function or the hyperparameters; see~\cite{jrnl_wang}. 
\end{remark}

\begin{algorithm}[t!]
\textbf{Initialize}: $\bx^1$, $\blambda^1 = \mathbf{0}$, $\nu^1 = 0$,  set the constant step size $\alpha$; prior on $\{\hat{U}^1_{m,n}\}$ if available. \\
\nl \bf for $t = 1, 2, \dots $ \bf do \\
\nl \quad \normalfont{Collect measurement of} ${\by}^t$ \\ 
\nl \quad \bf for \normalfont{each system $m = 1, 2, \dots M$} \bf do \\
\nl \quad \quad \normalfont{Update the input $x_{m}^t$ via~\eqref{eq:OPD_sys}} \\
\nl \quad \quad \normalfont{Send $x_{m}^t$ to users $n \in \mathcal{N}_m$} \\
\nl \quad \quad \normalfont{Update $\nu^t$ via~\eqref{eq:OPD_nu}} \\
\nl \quad \quad \bf for \normalfont{each user $n \in \mathcal{N}_m$} \bf do \\ 
\nl \quad \quad \quad \bf if \normalfont{Feedback is given by user $n$ } \\ 
\nl \quad \quad \quad \quad \normalfont{Update} $\hat{U}_{m,n}^t$ via \eqref{eq:GP_def_mu} \\
\nl \quad \quad \quad \bf else  \\
\nl \quad \quad \quad \quad \normalfont{Keep} $\hat{U}_{m,n}^t = \hat{U}_{m,n}^{t-1}$ \\
\nl \quad \quad \quad {\bf end if}  \\
\nl \quad \quad \quad \normalfont{Estimate} $g_{m,n}^t$ via \eqref{eq:gradient_U} \\ 
\nl \quad \quad \quad \normalfont{Update $x_{m,n}^t$ via \eqref{eq:OPD_user}} \\
\nl \quad \quad \quad \normalfont{Update $\lambda_{m,n}^t$ via~\eqref{eq:OPD_lambda}}\\
\nl \quad \quad \quad \normalfont{Send $x_{m,n}^t$ to system $m$ } \\
\nl \quad \quad  {\bf end for} (users)  \\
\nl \quad  \normalfont{System $m$ updates $\lambda_{m,n}^t$ via~\eqref{eq:OPD_lambda}}\\
\nl \quad  {\bf end for} (systems)  \\
\nl {\bf end for} (time)  \\
    \caption{{GP-based online primal-dual method} 
    \label{algo}}
\end{algorithm}

\subsection{Distributed GP-based Online Algorithm} 
\label{sec:algorithm}

\begin{figure}[t!] 
    \centering
\includegraphics[scale=0.48]{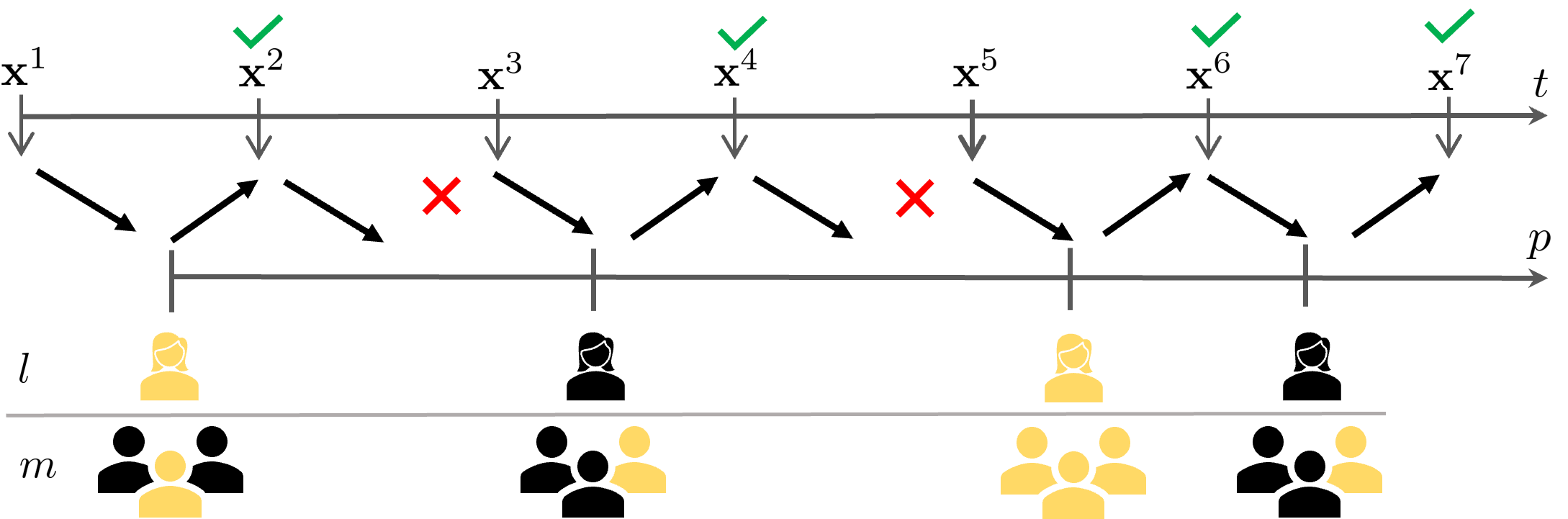}
    \vspace{0.1in}
   \caption{Feedback (color-coded in yellow) from users may come sporadically, and each user may provide feedback at different times (green ticks indicate that a functional evaluation  is received, red crosses indicate that a functional evaluation is not received). The functions $\{U_{m,n}\}$ are estimated concurrently with the execution of the Algorithm~\ref{algo}.} 
    \label{fig:algo_ex}
\end{figure}

Using the shape-constrained GPs machinery described in the previous section, the idea is to utilize the posterior mean~\eqref{eq:GP_def_mu} as a surrogate for $U_{m,n}$. In particular, for a given user $n \in \mathcal{N}_m$, the posterior mean~\eqref{eq:GP_def_mu} is computed based on the $p_{m,n}^t$ functional evaluations $\{U_{m,n}(x_{m}^{t_i}) + \epsilon_i\}_{i = 1}^{p_{m,n}^t}$, $t_i \in \{1, \dots, t\}$,  provided by the user   up to time $t$; that is, we set $\hat{U}^t_{m,n}(x) = \bar{\mu}_{p_{m,n}^t}({x})$. As pictorially illustrated in Figure~\ref{fig:algo_ex}, functional evaluations are provided sporadically by each user,  and each user may provide feedback at different times. Furthermore, feedback is not required at each step of the online algorithm as in zeroth-order methods.  

Once the posterior mean is available, an estimate of the gradient can  be obtained via, e.g., multivariable zeroth-order methods when $\bx^t_{m} \in \bbR^n$, or via finite difference when $x^t_{m} \in \mathbb{R}$. As an example for the latter, one way to obtain $g_{m,n}^t $ is: 
\vspace{-0.15cm}
\begin{equation}
\label{eq:gradient_U}
    g_{m,n}^t = \frac{\hat{U}^t_{m,n}(x^t_{m} + \delta/2) - \hat{U}^t_{m,n}(x^t_{m}  - \delta/2)}{\delta} 
\end{equation}
with $\delta$ a preselected parameter. 

The overall algorithm with feedback from both, users and network, is tabulated as Algorithm~\ref{algo}. In terms of operation, the updates~\eqref{eq:OPD_sys} and~\eqref{eq:OPD_nu} are implemented at each system $m$, with $\{x_m^t\}$ being physical inputs for the system; on the other hand, the update~\eqref{eq:OPD_user} is implemented at the user's side. The updates~\eqref{eq:OPD_lambda} are implemented by both users and systems, and users and systems are required to exchange the local variables $x_m^t$ and $x_{m,n}^t$.

\section{Performance Analysis}
\label{ref:performance}

In this section, we investigate the convergence of the online algorithm presented in Section \ref{sec:distri_PD}. To this end, the following standard assumptions are presumed. 

\begin{assumption}
\label{as:set_x} 
The set $\calX^t$ is convex and compact for all $t$.
\end{assumption}
\vspace{0.02in}
\begin{assumption}
\label{as:fun_f} 
The functions $\bx \mapsto f(\bx)$ and $\bx \mapsto C^t(\bar{\bA}^t \bx + \bB^t \bw^t)$ are convex and continuously differentiable for all $\{x_m, x_{m,n} \in \calX_m^t\}$. 
\end{assumption}
\vspace{0.02in}

The previous assumptions imply that $f$ and $C^t$ are Lipschitz continuous, their gradients are also bounded by their Lipschitz constant, and they are uniformly bounded for all $\bx^t \in \calX^t$ and for all $t$. Let $L$ be the Lipschitz constant of $f$. The next assumption pertains to each user-system star network. 
\begin{assumption}
\label{as:net} 
Each system $m$ is connected with its users via a star network (induced by~\eqref{eq:P1_c2}). The largest singular value of the incidence matrix $\bD$ is $\Omega$. 
\end{assumption} 
\vspace{0.02in}

Regarding Assumption~\ref{as:net}, each system-user network has a  diameter $h_m \leq 2$.

\vspace{0.02in}
\begin{assumption}
\label{as:slater} 
Slater’s constraint qualification holds $\forall~t$. 
\end{assumption} 
\vspace{0.02in}

\begin{assumption}
\label{as:set_la} 
(\cite{Online_PD}) The convex set $\Lambda^t$ is included in a 2-norm ball of radius $B_{\lambda}$, where $B_\lambda \geq N_{\max} h_{\max} L M + 1$, $h_{\max} = \underset{i=1,...,M}{\max}\{h_i\}$, and $N_{\max} =\underset{i=1,...,M}{\max}\{N_i\}$ for all $t$. 
\end{assumption}
\vspace{0.02in}
\begin{assumption}
\label{as:set_nu} 
(\cite{koshal2011multiuser}) The set $\Psi^t$ for the dual variable $\nu$ is convex and compact for all $t$.
\end{assumption}
\smallskip%

Under the Slater's constraint qualification, it is possible to compute the compact set $\Psi^t$ as shown in \cite{koshal2011multiuser}. On the other hand, an analytical expression for the radius of  $\Lambda^t$ is not available.  In a practical implementation of the algorithm, in
case of persistent and large consensus errors (which are due to a clipping of the dual variables if the radius of  $\Lambda^t$ is too small), the radius of  $\Lambda^t$ can be increased.

We tackle the analysis by viewing our  algorithm as an online primal-dual with \emph{inexact} gradient information~\cite{bertsekas2000gradient}. In particular, the approximate gradient (obtained based on a finite difference method based on the estimated functions $\{\hat{U}^t_{m,n}\}$) can be expressed as $\bg^t =  \nabla f(\bx^t) + \be^t$, where $\nabla f(\bx^t)$ is the true gradient and $\be^t$ is a stochastic vector. Before stating appropriate assumptions, define the filtration $\mathcal{F}_t = \{\be^1, \ldots, \be^{t-1}\}$. The following assumption is made. 


\vspace{.1cm}
\begin{assumption}
\label{as:error} 
$\exists~\bar{E}^t < \infty$ such that $\bbE[\|\be^t\|^2|\mathcal{F}_t] \leq \bar{E}^t$, $\forall~t$.
\end{assumption}  
\vspace{.1cm}



Notice that the error $\be^t$ is not assumed to have zero mean, since the GP may introduce a bias in the estimated gradient.  
Assumption~\ref{as:error} presumes that the expected value of $\|\be^t\|^2$ is bounded; this is consistent with standard assumptions in the context of inexact gradient methods and stochastic gradient methods (see, e.g.,~\cite{bertsekas2000gradient,gower2019sgd}). We also note that $\bar{E}^t$ can be taken to be the Bayes errors for regression
with GPs; in fact,~\cite{jrnl_wang} showed that the error of the GP regression is higher than the error incurred by the shape-constrained GP regression. We  point out that Assumption~\ref{as:error} is also supported by empirical evidences; an example for the shape-constrained GP is provided in Figure~\ref{fig:grad_error}, which illustrates the estimated function and the estimated derivative for different numbers of noisy functional evaluations. With these standard assumptions in place, the next section will provide performance bounds for the proposed algorithm.

\subsection{Dynamic Regret}

We characterize the performance of Algorithm~\ref{algo} using the so-called dynamic regret (see, e.g.,~\cite{Shahrampour2018,yi2016tracking,Online_PD,lee2016distributed} and references therein). In particular, given the distributed nature of our algorithm, we consider an extension of the network regret~\cite{Online_PD,lee2016distributed} to a time-varying setting. To this end, define the \textit{dynamic regret} per user $j \in \mathcal{N}_m$ interacting with the system $m$, with respect to an optimal solution of \eqref{eq:P1} for the system $m$ at time $t$, i.e., $x_m^{t*}$, as 
\begin{equation}
    \label{eq:reg_user}
    \mathrm{Reg}_T^{m,j} := \sum_{t=1}^T \sum_{i=1}^{N_m} {U}_{m,i}(x_{m,j}^t) - \sum_{t=1}^T \sum_{i=1}^{N_m} {U}_{m,i}(x_{m}^{t*}).
\end{equation}
Then, based on \eqref{eq:reg_user}, we define the \textit{dynamic network regret} for the system $m$ as~\cite{Online_PD,lee2016distributed}
\begin{align}
    \label{eq:reg_local} \nonumber
    & \mathrm{Reg}_T^{m} = \frac{1}{N_m} \sum_{j=1}^{N_m} \mathrm{Reg}_T^{m,j} \\ 
    & ~~~~ = \frac{1}{N_m} \sum_{t=1}^T \sum_{i,j=1}^{N_m} {U}_{m,i}(x_{m,j}^t) - \sum_{t=1}^T \sum_{i=1}^{N_m} {U}_{m,i}(x_{m}^{t*}).
\end{align}
By using \eqref{eq:reg_local}, the dynamic \textit{global network regret} after $T$ iterations can be defined as 
\begin{align}
    \label{eq:reg}
    \mathrm{Reg}_T & := \sum_{m=1}^{M} \mathrm{Reg}_T^{m}.
\end{align}
To analyze the global network regret, let $\bx^{t*}$ be an optimal solution of \eqref{eq:P1} at time $t$, and consider the following standard definitions of the path length, which captures the drift of the optimal solutions over $T$ time steps
\begin{align}
    \label{eq:def1}
    \hspace{-.2cm} \Phi^T := \sum_{t=1}^T \|\bx^{t*}-\bx^{t+1*}\| , \,\, \Upsilon^T := \sum_{t=1}^T \|\bx^{t*}-\bx^{t+1*}\|^2 . 
\end{align}
Further, consider the following definitions for the cumulative expected error in the gradient: 
\begin{align}
    \label{eq:def2}
    \xi^T := \sum_{t=1}^T \bbE[\|\be^t\| |\mathcal{F}_t] , \qquad \Xi^T := \sum_{t=1}^T \bar{E}^t.
\end{align}

With these metrics in place, we start by bounding the norm of the true gradient of the time-varying Lagrangian with respect to the primal variable for any $\bx \in \calX^t$, $\nu \in \Psi^t$ and $\blambda \in \Lambda^t$. Let $ J := \underset{t \in \calT}{\sup} \, \underset{\bx \in \calX^t}{\max} \| \nabla C^t(\by^t(\bx)) \|$ by Assumption \ref{as:set_x} and Assumption \ref{as:fun_f}, $B_\lambda := \underset{t \in \calT}{\sup} \, \underset{\lambda \in \Lambda^t}{\max} \|\blambda\|$ and $B_{\nu}:= \underset{t \in \calT}{\sup} \, \underset{\nu \in \Psi^t}{\max} \, |\nu |$ under Assumption \ref{as:set_la} and Assumption \ref{as:set_nu}, then we have
\begin{align} \nonumber 
    \|\nabla_{\bx} \calL^t(\bx , \nu, \blambda) \| &\leq \| \nabla f(\bx) \| + \nu \| \nabla C^t(\by^t(\bx))\| + \| \bD^\top \blambda \| \\ 
    &\leq L +  B_{\nu} J + \Omega B_{\lambda} := \Gamma_x,
    \label{eq:grad_bound_lagx}
\end{align}
where, the last inequality holds by using the triangle and Cauchy-Schwarz inequality, and Assumption \ref{as:net}. The term $\Gamma_x = L +  B_{\nu} J + \Omega B_{\lambda}$ is an upper bound for the gradient of the Lagrangian with respect to $\bx$.

Similarly, we bound the norm of the gradient of the Lagrangian with respect to the dual variables $\blambda$ and $\nu$ for any $\bx \in \calX^t$, $\nu \in \Psi^t$ and $\blambda \in \Lambda^t$. Let $\bkappa = [\nu, \blambda^\top]^\top$, $\left|C^t(\by^t(\bx)) \right| \leq H < \infty$ for all $t$ by Assumption \ref{as:set_x} and Assumption \ref{as:fun_f}, and $B_x := \underset{t \in \calT}{\sup} \, \underset{\bx \in \calX^t}{\max} \|\bx\|$ under Assumption \ref{as:set_x}. Then, $\|\nabla_{\kappa} \calL^t(\bx, \bkappa)\|^2$ can be bounded by
\begin{align} \nonumber
    \|\nabla_{\kappa} \calL^t(\bx, \bkappa)\|^2 &= \left\| \begin{matrix}
    \nabla_{\lambda} \calL^{t}(\bx , \nu, \blambda) \\ \nabla_{\nu} \calL^{t}(\bx, \nu, \blambda)
    \end{matrix} \right\|^2, \\
    &\leq \Omega^2 B_x^2 + H^2 := \Gamma_\kappa, 
    \label{eq:P9_T1_aa}
\end{align}
where we used the Cauchy-Schwarz inequality, Assumption \ref{as:set_x}, and Assumption \ref{as:net}. Further, let $\mathrm{diam}(\calX^t) \leq D^t_x$ where $D^t_x \leq D_x < \infty$ by Assumption \ref{as:set_x}. The term $\Gamma_\kappa = \Omega^2 B_x^2 + H^2$ is an upper bound for the gradient of the Lagrangian with respect to $\bkappa$.

With these notations in place, we are now ready to state the main result for the dynamic regret. 

\vspace{.2cm}

\begin{theorem} \label{th1}
Let Assumptions \ref{as:set_x}-\ref{as:error} hold. Set $\blambda^1 = \mathbf{0}$, $\nu^1 = 0$. Then, the dynamic network regret after $T$ iterations is upper bounded by
\begin{align} \nonumber
    &\bbE[\mathrm{Reg}_T] \leq \frac{1}{2 \alpha} (\| \bx^1 - \bx^{1*} \|^2 + B_\lambda^2 + B_\nu^2) + \frac{\alpha}{2} T (\Gamma_x^2 + \Gamma_\kappa)\\ 
    &\quad + \frac{\alpha}{2} \Xi^T + \xi^T \, (2B_x + \alpha \Gamma_x) + \frac{1}{2 \alpha} \Upsilon^T + \frac{1}{\alpha} D_x \Phi^T.
    \label{eq:P17_T1}
\end{align}
\end{theorem}

\vspace{.2cm}

The proof of this theorem will be provided in Section~\ref{sec:proofs}. Theorem \ref{th1} asserts that the bound on the expected dynamic network regret $\bbE[\mathrm{Reg}_T]$ after $T$ time steps depends on the temporal variability of the solutions (through the terms $\frac{1}{2 \alpha} \Upsilon^T$ and $\frac{1}{\alpha} D_x \Phi^T$) and the errors associated with the estimation of the gradient (through the terms $\frac{\alpha}{2} \Xi^T$ and $\xi^T (2B_x + \alpha \Gamma_x)$). The remaining terms are in line with standard regret results (e.g.,~\cite{Online_PD}). In particular, we can observe the following: 

\vspace{.1cm}

\noindent $\bullet$ \emph{Persistent variations and gradient errors:} Suppose that $\Upsilon^T$ and $\Phi^T$ grow as $\mathcal{O}(T)$, modeling a persistent temporal variation of the constraints and/or the vector $\bw^t$; suppose further that $\Xi^T$ and $\xi^T$ grow as $\mathcal{O}(T)$, modeling a persistent error in the gradient $\bg^t$. Then, $\bbE[\mathrm{Reg}_T]/T$ behaves  as $\mathcal{O}(1)$.    
\vspace{.1cm}

\noindent  $\bullet$ \emph{Vanishing gradient errors:} If $\Xi^T$ and $\xi^T$ grows sublinearly, i.e., as $o(T)$, then the asymptotic average regret is in general 
$\frac{1}{T}\bbE[\mathrm{Reg}_T] = \mathcal{O}(1 + T^{-1} \Upsilon^T + T^{-1} \Phi^T)$.
Notice that if the algorithm is executed over a finite interval of $T$ steps and the step-size is chosen as $\alpha = \frac{1}{\sqrt{T}}$, then $\bbE[\mathrm{Reg}_T]/T$ behaves as 
$\frac{1}{T}\bbE[\mathrm{Reg}_T] = \mathcal{O}(T^{-\frac{1}{2}} + T^{-\frac{1}{2}} \Upsilon^T + T^{-\frac{1}{2}} \Phi^T)$ \, .

\vspace{.1cm}

\noindent  $\bullet$ \emph{Vanishing variations and errors:} If $\Upsilon^T$, $\Phi^T$, $\Xi^T$, and $\xi^T$ all grow as $o(T)$, then we recover the asymptotic result of~\cite{Online_PD}.  

\vspace{.1cm}

\begin{figure}[t] 
    \centering
\includegraphics[scale=0.525]{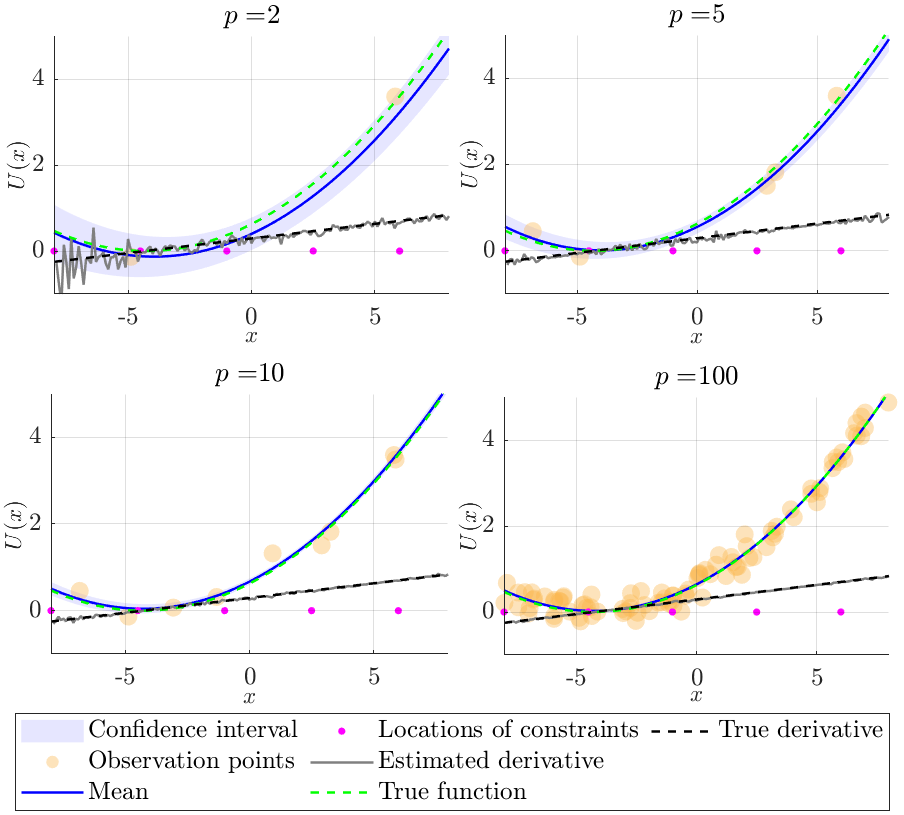}
   \caption{Example of the derivative estimation using shape-constrained GPs with hyperparameters $\sigma_f = 1$ and $\ell=10$. As expected, the estimated derivative error decreases with the increasing of the number of observations $p$. In this case, shape-constrained GP outperforms GP for low numbers of samples (see also Section~\ref{sec:results}).}
    \label{fig:grad_error}
    \vspace{-.4cm}
\end{figure}

\begin{remark}
Regarding the errors in the gradient, if the function $U_{m,n}$ is strongly convex, then the estimation error yielded by a shape-constrained GPs will decrease with the increasing of the number of functional evaluations; see the illustrative example in Figure~\ref{fig:grad_error}, and the discussion in~\cite{jrnl_wang}. However, zeroth-order and finite-different methods will still generate an error in the estimation of the gradient, and thus $\Xi^T$ and $\xi^T$  would generally increase as $\calO(T)$. To obtain a trend as  $o(T)$, one has to increase the accuracy in the gradient estimation to make $\be^t$ arbitrarily small. \hfill $\Box$
\end{remark}

\subsection{Average Constraint Violation}

We now consider a bound on the violation of the constraint on the network output $\by^t$. To this end, we consider the so-called average constraint violation (ACV), which we define here as \vspace{-0.25cm}
\begin{align}
    \label{eq:fit}
    \mathrm{ACV}_T := \sum_{t=1}^T \left[ C^{t}(\by^{t}(\bx^t)) \right]^+.
\end{align}
We note that some prior works in context (e.g.,~\cite{chen2018bandit,yi2020distributed}) consider the modified definition $\left[ \sum_{t=1}^T C^{t}(\by^{t}(\bx^t)) \right]^+$, which may lead to looser bounds. 

\vspace{.2cm}

\begin{theorem} \label{th2}
Let Assumptions \ref{as:set_x}-\ref{as:error} hold. Then, the average constraint violation incurred by the algorithm \eqref{eq:OPD_case1} can be bounded as
\begin{align} \nonumber
    &\bbE \left[ \mathrm{ACV}_T \right] \leq B_\nu^{-1} T \left( D_x L + B_\lambda \Omega B_x \right) \\ \nonumber
    & ~~~ + B_\nu^{-1} \left( \frac{\alpha}{2} \Xi^T + \xi^T (2B_x + \alpha \Gamma_x) + \frac{1}{2 \alpha} \Upsilon^T + \frac{1}{\alpha} D_x \Phi^T \right) \\
    & ~~~ + B_\nu^{-1} T \left( \frac{1}{\alpha} (4 B_x^2 +  B_\nu^2) + \frac{\alpha}{2} (\Gamma_x^2 + H^2) \right).
    \label{eq:P5_T2a}
\end{align}
\end{theorem}

\vspace{.2cm}

The second and third terms of \eqref{eq:P5_T2a} exhibit the same asymptotic behavior of the dynamic regret under persistent variations and gradient errors, and with vanishing variations and errors. In case of vanishing variations and gradient errors, the bottlenecks are the two terms $T D_x L$ and $T B_\lambda \Omega B_x$, which make $\bbE[\mathrm{ACV}_T]/T$ behave as $\mathcal{O}(1)$.   


\subsection{Proofs of the Results}
\label{sec:proofs}

To streamline exposition,  define $\phi^t := \|\bx^{t*}-\bx^{t+1*}\|$ and $\bar{\mathrm{e}}^t := \bbE[\|\be^t\|| \mathcal{F}_t]$. To derive our main results, we will utilize the following lemma.

\vspace{.2cm}

\begin{lemma} \label{lemma1}
Let Assumptions \ref{as:set_x}-\ref{as:error} hold, and set $\blambda^1 = \mathbf{0}$, $\nu^1 = 0$. Then, for any $\bkappa \in \Psi^t \times \Lambda^t$ the following holds:
\begin{align} \nonumber
    &\bbE\left[\sum_{t=1}^T  \left( \calL^t(\bx^t,  \bkappa) - \calL^t(\bx^{t*} , \bkappa^t) \right) \right]\\ \nonumber
    &\quad \leq  \frac{1}{2 \alpha} (\| \bx^1 - \bx^{1*} \|^2 + \|\bkappa\|^2) + \frac{\alpha}{2} T (\Gamma_x^2 + \Gamma_\kappa) \\ 
    &\quad + \frac{\alpha}{2} \Xi^T  + \xi^T \, (2 B_x + \alpha \Gamma_x) + \frac{1}{2 \alpha} \Upsilon^T + \frac{1}{\alpha} D_x \Phi^T.
    \label{eq:P10b_T1}
\end{align}
\end{lemma}

\vspace{.2cm}

\textit{Proof.} Notice first that
\begin{align} 
    \|\tilde \nabla_{\bx} \calL^t(\bx^t , \nu, \blambda) - \nabla_{\bx} \calL^t(\bx^t, \nu , \blambda) \|  &= \| \be^t \|,
    \label{eq:grad_bound_difflagx} 
\end{align}
where $\tilde \nabla_{\bx} \calL^t(\bx^t , \nu, \blambda)$ is the inexact gradient of the Lagrangian function (where we utilize $\bg^t$ instead of $\nabla f(\bx^t)$). Recall that $\bx^{t *}$ denotes an optimal solution at time $t$.  Using the primal update, adding and subtracting $\bx^{t *}$, we get: 
\begin{align} \nonumber
    &\|\bx^{t+1} - \bx^{t+1*}\|^2 = \|\bx^{t+1} - \bx^{t*} + \bx^{t*} - \bx^{t+1*}\|^2 \\ \nonumber
    &~~ = \|\bx^{t+1} - \bx^{t*} \|^2 + 2(\bx^{t+1} - \bx^{t*})^\top (\bx^{t*} - \bx^{t+1*}) \\ \nonumber
    &~~  \quad + \|\bx^{t*} - \bx^{t+1*}\|^2, \\ \nonumber
    &~~  = \left\| \mathrm{proj}_{{\calX^t}} \left\{ \bx^t -\alpha \tilde \nabla_{\bx} \calL^t(\bx^t , \nu^t, \blambda^t) \right\} - \bx^{t*} \right\|^2 + (\phi^t)^2 \\ \nonumber
    &~~  \quad + 2(\bx^{t+1} - \bx^{(*,t)})^\top (\bx^{t*} - \bx^{t+1*}), \\ \nonumber
    &~~  \leq \|  \bx^t -\alpha \tilde \nabla_{\bx} \calL^t(\bx^t, \nu^t, \blambda^t) - \bx^{t*} \|^2 + (\phi^t)^2 \\ \nonumber 
    &~~  \quad + 2 \|\bx^{t+1} - \bx^{t*}\| \|\bx^{t*} - \bx^{t+1*}\|, \\ 
    &~~  \leq \| \bx^t - \bx^{t*} -\alpha \tilde \nabla_{\bx} \calL^t(\bx^t, \nu^t, \blambda^t) \|^2 + (\phi^t)^2 + 2 D_x \phi^t,
    \label{eq:P1_T1}
\end{align}
where \eqref{eq:P1_T1} holds by the non-expansiveness property of the projection operator, the Cauchy-Schwarz inequality, and the definition of $\phi^t$.

By adding and subtracting $\alpha \nabla_\bx \calL^t(\bx^t, \nu^t, \blambda^t)$  (to simplify notation we use $\nabla_\bx \calL^t=\nabla_\bx \calL^t(\bx^t, \nu^t, \blambda^t)$), and expanding the first term of the right-hand side (RHS) of \eqref{eq:P1_T1}, we get
\begin{align} \nonumber
    &\| \bx^t - \bx^{t*} -\alpha \tilde \nabla_{\bx} \calL^t \|^2  \\ \nonumber
    &\quad = \| \bx^t - \bx^{t*} - \alpha \nabla_{\bx} \calL^t + \alpha (\nabla_{\bx} \calL^t - \tilde \nabla_{\bx} \calL^t)\|^2, \\ \nonumber
    &\quad= \| \bx^t - \bx^{t*} - \alpha \nabla_{\bx} \calL^t \|^2 + \alpha^2 \|\nabla_{\bx} \calL^t - \tilde \nabla_{\bx} \calL^t \|^2 \\ \nonumber
    &\quad \quad + 2\alpha (\bx^t - \bx^{t*} - \alpha \nabla_{\bx} \calL^t)^\top  (\nabla_{\bx} \calL^t - \tilde \nabla_{\bx} \calL^t ),\\ \nonumber
    &\quad \leq \| \bx^t - \bx^{t*} - \alpha \nabla_{\bx} \calL^t \|^2 + \alpha^2 \|\nabla_{\bx} \calL^t - \tilde \nabla_{\bx} \calL^t\|^2 \\
    &\quad \quad + 2 \alpha \|\bx^t - \bx^{t*} - \alpha \nabla_{\bx} \calL^t\|  \| \nabla_{\bx} \calL^t - \tilde \nabla_{\bx} \calL^t\|.
    \label{eq:P2_T1_a}
\end{align}

By using \eqref{eq:grad_bound_difflagx},  we have that \eqref{eq:P2_T1_a} can be bounded as
\begin{align} \nonumber
     &\| \bx^t - \bx^{t*} -\alpha \tilde \nabla_{\bx} \calL^t \|^2 \\ \nonumber
     &\quad \leq \| \bx^t - \bx^{t*} - \alpha \nabla_{\bx} \calL^t \|^2 + \alpha^2 \| \be^t \|^2 \\ \nonumber
     &\quad \quad + 2 \alpha \| \be^t \| \, \|\bx^t - \bx^{t*} - \alpha \nabla_{\bx} \calL^t\|, \\ \nonumber
     &\quad \leq \| \bx^t - \bx^{t*} - \alpha \nabla_{\bx} \calL^t\|^2 + \alpha^2 \| \be^t \|^2 \\\nonumber
     &\quad \quad + 2 \alpha \| \be^t \| \, (\|\bx^t - \bx^{t*} \| + \alpha \| \nabla_{\bx} \calL^t\|), \\  
     &\quad \leq \| \bx^t - \bx^{t*} - \alpha \nabla_{\bx} \calL^t \|^2 + \alpha^2 \| \be^t \|^2 + 2 \alpha \| \be^t \| \, (2B_x + \alpha \Gamma_x),
    \label{eq:P2_T1}
\end{align}
where the last inequality is derived using \eqref{eq:grad_bound_lagx}  and the triangle inequality. By expanding the first term on the RHS of \eqref{eq:P2_T1}, we get
\begin{align} \nonumber
    &\| \bx^t - \bx^{t*} - \alpha \nabla_{\bx} \calL^t \|^2 \\ \nonumber
    &\quad = \| \bx^t - \bx^{t*} \|^2 + \alpha^2 \|\nabla_{\bx} \calL^t\|^2 - 2 \alpha (\bx^t - \bx^{t*})^\top \nabla_{\bx} \calL^t, \\ \nonumber
    &\quad \leq \| \bx^t - \bx^{t*} \|^2 + \alpha^2 \Gamma_x^2 \\
    &\quad \quad - 2 \alpha [ \calL^t(\bx^t, \nu^t, \blambda^t)  - \calL^t(\bx^{t*} , \nu^t, \blambda^t)].
    \label{eq:P4_T1}
\end{align}
Since the function $f(\bx)$ is convex, the time-varying Lagrangian is convex in $\bx$. Thus, from the first-order characterization of convexity \eqref{eq:P4_T1} holds.

Using \eqref{eq:P4_T1} the bound  \eqref{eq:P2_T1}, one has that 
\begin{align} \nonumber
   & \| \bx^t - \bx^{t*} -\alpha \tilde \nabla_{\bx} \calL^t\|^2 \\ \nonumber
   &\quad \leq \| \bx^t - \bx^{t*} \|^2 + \alpha^2 \Gamma_x^2 + \alpha^2 \| \be^t \|^2 + 2 \alpha \| \be^t \| \, (2B_x + \alpha \Gamma_x) \\
   &\quad \quad  - 2 \alpha [ \calL^t(\bx^t, \nu^t, \blambda^t)  - \calL^t(\bx^{t*} , \nu^t, \blambda^t) ].
    \label{eq:P5_T1}
\end{align}
By using \eqref{eq:P5_T1} in \eqref{eq:P1_T1}, we have
\begin{align} \nonumber
    &\|\bx^{t+1} - \bx^{t+1*}\|^2 \\\nonumber
    & \leq \| \bx^t - \bx^{t*} \|^2 + \alpha^2 \Gamma_x^2 + \alpha^2 \| \be^t \|^2 + 2 \alpha \| \be^t \| \, (2B_x + \alpha \Gamma_x) \\
    &  - 2 \alpha [ \calL^t(\bx^t, \nu^t, \blambda^t)  - \calL^t(\bx^{t*} , \nu^t, \blambda^t) ] + (\phi^t)^2 + 2 D_x \phi^t .
    \label{eq:P6_T1}
\end{align}
It thus follows that: 
\begin{align} \nonumber
    &\calL^t(\bx^t, \nu^t, \blambda^t)  - \calL^t(\bx^{t*} , \nu^t, \blambda^t) \\\nonumber 
    & \leq \frac{1}{2 \alpha} (\| \bx^t - \bx^{t*} \|^2 - \|\bx^{t+1} - \bx^{t+1*}\|^2) + \frac{\alpha}{2} \Gamma_x^2 \\
    & + \frac{\alpha}{2} \| \be^t \|^2 + \| \be^t \| \, (2B_x + \alpha \Gamma_x) + \frac{1}{2 \alpha} (\phi^t)^2 + \frac{1}{\alpha} D_x \phi^{t}.
    \label{eq:P7_T1}
\end{align}
Now, we proceed with a similar analysis for the distance between the updates of the dual variables $\blambda^{t+1}, \nu^{t+1}$ and an arbitrary dual variable, i.e., for any $\blambda \in \Lambda^t $ and any $ \nu \in \Psi^t$. We start with the following: 
\begin{align} \nonumber
    &\|\bkappa^{t+1} - \bkappa\|^2 \leq \|\bkappa^t + \alpha \nabla_{\kappa} \calL^t - \bkappa\|^2, \\ \nonumber
    &\quad \leq \|\bkappa^t - \bkappa\|^2 + \alpha^2 \|\nabla_{\kappa} \calL^t\|^2 + 2 \alpha \nabla_{\kappa} \calL^{t\top} (\bkappa^t - \bkappa), \\
    &\quad \leq \|\bkappa^t - \bkappa\|^2 + \alpha^2 \Gamma_\kappa + 2 \alpha \nabla_{\kappa} \calL^{t\top} (\bkappa^t - \bkappa) 
    \label{eq:P8_T1_aa}
\end{align}
where $\nabla_{\kappa} \calL^t = \nabla_{\kappa} \calL^t(\bx^t, \bkappa^t)$ for brevity. 

The time-varying Lagrangian is concave in $\blambda$ and $\nu$; i.e., any $\bkappa$, which implies that the time-varying Lagrangian differences for a fixed $\bx^t \in \calX^t$ satisfies
\begin{align}
    \calL^t(\bx^t, \bkappa^t)  - \calL^t(\bx^t, \bkappa) &\geq  (\bkappa^t - \bkappa)^\top \nabla_{\kappa} \calL^t(\bx^t, \bkappa).
    \label{eq:lag_lin_aa}
\end{align}
Substituting \eqref{eq:lag_lin_aa} in \eqref{eq:P8_T1_aa} for any $\bkappa \in \Psi^t \times \Lambda^t$, we get
\begin{align} \nonumber
    &\calL^t(\bx^t, \bkappa^t)  - \calL^t(\bx^t, \bkappa) \\
    &\quad \geq \frac{1}{2 \alpha} (\|\bkappa^{t+1} - \bkappa\|^2 - \|\bkappa^t - \bkappa\|^2) - \frac{\alpha}{2} \Gamma_\kappa.
    \label{eq:P10_T1_aa}
\end{align}
By subtracting the results in \eqref{eq:P7_T1} and \eqref{eq:P10_T1_aa} we have
\begin{align} \nonumber
    &\calL^t(\bx^t, \bkappa) - \calL^t(\bx^{t*}, \bkappa^t) \\\nonumber
    &\quad \leq \frac{1}{2 \alpha} (\| \bx^t - \bx^{t*} \|^2 - \|\bx^{t+1} - \bx^{t+1*}\|^2) \\\nonumber
    &\quad \quad +  \frac{1}{2 \alpha} (\|\bkappa^t - \bkappa\|^2-\|\bkappa^{t+1} - \bkappa\|^2) + \| \be^t \| \, (2B_x + \alpha \Gamma_x) \\ 
    &\quad \quad + \frac{\alpha}{2} \| \be^t \|^2 + \frac{1}{2 \alpha} (\phi^t)^2 + \frac{1}{\alpha} D_x \phi^t + \frac{\alpha}{2} (\Gamma_x^2 + \Gamma_\kappa).
    \label{eq:P10_T1}
\end{align}
Assume that  $\bkappa^{1} = \mathbf{0}$. Then, summing \eqref{eq:P10_T1} over time, using the telescopic property for the first two terms on the RHS, and by taking the expectation of both sides, the result follows. 
\hfill $\blacksquare$ 

Based on Lemma \ref{lemma1}, in the following, we provide the proof for Theorem~\ref{th1}. 

\vspace{.2cm}

\subsubsection{Proof for Theorem~\ref{th1}}

First, we can rewrite $\calL^t(\bx^t, \bkappa) - \calL^t(\bx^{t*} , \bkappa^t)$ as:
\begin{align} \nonumber
    &\calL^t(\bx^t, \bkappa) - \calL^t(\bx^{t*} , \bkappa^t) = f(\bx^t) - f(\bx^{t*}) \\
    & ~~~~~~~~~ + \underbrace{ \bkappa^\top \left[ \begin{matrix} C^t \left( \by^t(\bx^t) \right) \\ \bD \bx^t \end{matrix} \right]  - \bkappa^{t \top} \left[ \begin{matrix} C^t \left( \by^t(\bx^{t*}) \right) \\ \bD \bx^{t*} \end{matrix} \right]}_{:= \mathrm{\beta}^t}
    \label{eq:P11_T1}
\end{align}
where we recall that $\bD \bx^{t*} = \mathbf{0}$ and $C^t \left( \by^t(\bx^{t*}) \right) \leq 0$. 

By using the definition of the cost function $f(\bx)$ and summing \eqref{eq:P11_T1} over time, we have
\begin{align} \nonumber
    &\sum_{t=1}^T \left( \calL^t(\bx^t, \bkappa) - \calL^t(\bx^{t*} , \bkappa^t) \right) = \sum_{t=1}^T \beta^t \\
    & + \sum_{t=1}^T \left(  \sum_{m=1}^M \sum_{n=1}^{N_m} {U}_{m,n}(x_{m,n}^t) - \sum_{m=1}^M \sum_{n=1}^{N_m} {U}_{m,n}(x_{m}^{t*}) \right).
    \label{eq:P11_T1a}
\end{align}

Taking the expected value and using Lemma \ref{lemma1} in \eqref{eq:P11_T1a}, we have
\begin{align} \nonumber 
    &\hspace{-0.6cm} \bbE\left[ \sum_{t=1}^T \left(  \sum_{m=1}^M \sum_{n=1}^{N_m} {U}_{m,n}(x_{m,n}^t) - \sum_{m=1}^M \sum_{n=1}^{N_m} {U}_{m,n}(x_{m}^{t*}) \right) + \sum_{t=1}^T \beta^t \right] \\ \nonumber
    &\quad \leq \frac{1}{2 \alpha} (\| \bx^1 - \bx^{1*} \|^2 + \|\bkappa\|^2) + \frac{\alpha}{2} T (\Gamma_x^2 + \Gamma_\kappa)\\ 
    &\quad \quad + \frac{\alpha}{2} \Xi^T + \xi^T \, (2B_x + \alpha \Gamma_x) + \frac{1}{2 \alpha} \Upsilon^T + \frac{1}{\alpha} D_x \Phi^T.
    \label{eq:P11a_T1}
\end{align}

In order to bound the global network regret, we add and subtract the RHS of \eqref{eq:P11_T1a} in \eqref{eq:reg}, to obtain
\begin{align} \nonumber
    &\mathrm{Reg}_T =  \\ \nonumber
    &\sum_{t=1}^T \left[ \sum_{m=1}^{M} \left(  \frac{1}{N_m}  \sum_{i,j=1}^{N_m} {U}_{m,i}(x_{m,j}^t)  -  \sum_{n=1}^{N_m} {U}_{m,n}(x_{m,n}^t) \right) - \beta^t \right] \\ 
    &+ \sum_{t=1}^T \left[  \left( \sum_{m=1}^M \sum_{n=1}^{N_m} {U}_{m,n}(x_{m,n}^t) - \sum_{m=1}^M \sum_{n=1}^{N_m} {U}_{m,n}(x_{m}^{t*}) \right) + \beta^t \right]. 
    \label{eq:P12_T1}
\end{align}
The expected value of the last term of \eqref{eq:P12_T1} is bounded by \eqref{eq:P11a_T1}. Next, focus on the first term in \eqref{eq:P12_T1}. Assumption \ref{as:fun_f} implies that
\begin{align} \nonumber
    &\sum_{i,j=1}^{N_m} {U}_{m,i}(x_{m,j}^t)  -  \sum_{n=1}^{N_m} {U}_{m,n}(x_{m,n}^t) \\ 
    &\quad \leq \sum_{i,j=1}^{N_m} L_{m,i} |x_{m,j}^t - x_{m,i}^t|,
    \label{eq:P13_T1}
\end{align}
where $L_{m,i} \leq L_m \leq L$ per each network composed of a system and the users interacting with such a system. Maximizing over the RHS of \eqref{eq:P13_T1} we can get an expression for the worst-case disagreement per system $m$ as
\begin{align} 
    \sum_{i,j=1}^{N_m} L_{m,i} |x_{m,j}^t- x_{m,i}^t|  \leq N_m^2 L_m h_m \| \bD_m \bx_m^t\|,
    \label{eq:P14_T1}
\end{align}
where $\underset{i,j}{\max} |x_{m,j}^t - x_{m,i}^t|$ can by upper bounded using Assumption \ref{as:net}, and $\bx_m := [x_m, x_{m,1}, \dots, x_{m,N_m}] \, \forall \, m$.

Back to the first term in \eqref{eq:P12_T1}, we have that
\begin{align} \nonumber
    &\sum_{t=1}^T \left[ \sum_{m=1}^{M} \left(  \frac{1}{N_m}  \sum_{i,j=1}^{N_m} {U}_{m,i}(x_{m,j}^t)  -  \sum_{n=1}^{N_m} {U}_{m,n}(x_{m,n}^t) \right) - \beta^t \right] \\ \nonumber
    &\quad = \sum_{t=1}^T \left[ \sum_{m=1}^{M} \frac{1}{N_m} \sum_{i,j=1}^{N_m} \left( {U}_{m,i}(x_{m,j}^t) - {U}_{m,n}(x_{m,n}^t) \right) - \beta^t \right] \\ \nonumber
    &\quad \leq  \sum_{t=1}^T  \sum_{m=1}^{M} \frac{1}{N_m} (N_m^2 L_m h_m \| \bD_m \bx_m^t\| ) - \sum_{t=1}^T  \beta^t\\ 
    %
    &\quad  \leq  \sum_{t=1}^T \, N_{\max} h_{\max}  L \, M \|\bD \bx^t\| - \sum_{t=1}^T \beta^t \, .
    \label{eq:P15_T1}
\end{align}
By the fact that at optimality $\blambda^{t \top} \bD \bx^{t*} = 0$ and $\nu^t C^t (\by^t(\bx^{t*})) \leq 0$, we can rewrite the RHS of \eqref{eq:P15_T1} as
\begin{align} \nonumber
    &\sum_{t=1}^T N_{\max}  h_{\max}  L M \|\bD \bx^t\| - \sum_{t=1}^T \bkappa^\top \left[ \begin{matrix}  C^t \left( \by^t(\bx^t) \right) \\ \bD \bx^t \end{matrix} \right] \\ \nonumber
    &\quad =\sum_{t=1}^T \left( N_{\max}  h_{\max}  L M \frac{\bD \bx^t}{\|\bD \bx^t\|} - \blambda \right)^\top \bD \bx^t \\
    &\quad \quad - \sum_{t=1}^T \nu \, C^t \left( \by^t(\bx^t)\right).
    \label{eq:P16_T1}
\end{align}
By construction, we can choose a feasible dual variable $\bar{\blambda}$ using the compactness of $\calX^t$ and the boundedness of $\| \bD \bx^t \|$, as in \cite{Online_PD}, so that the first term in \eqref{eq:P16_T1} equals $0$. Similarly, we can set $\bar{\nu}=0$. With this choice for  $\bar{\bkappa} = [\bar{\nu}, \bar{\blambda}^\top]^\top$, the expression on \eqref{eq:P16_T1} is equal to $0$.  Then, the expected value of the  regret defined in \eqref{eq:reg} is bounded as shown in \eqref{eq:P17_T1}. \hfill $\blacksquare$

\vspace{.2cm}

\subsubsection{Proof for Theorem~\ref{th2}}


We start by bounding the distance between the update of the dual variable $\nu^{t+1}$ and an arbitrary dual variable, i.e., for any $\nu \in \Psi^t$ 
\begin{align} \nonumber
    |\nu^{t+1} - \nu|^2 &\leq |\nu^t + \alpha \nabla_{\nu} \calL^t - \nu|^2, \\ 
    &\leq |\nu^t - \nu|^2 + \alpha^2 H^2 + 2 \alpha \nabla_{\nu} \calL^{t} (\nu^t - \nu),
    \label{eq:P1_T2r}
\end{align}
where $\nabla_{\nu} \calL^t = \nabla_{\nu} \calL^t(\bx^t, \nu, \blambda) = C^t(\by^t(\bx^t))$ for any feasible $\blambda$. Reorganizing \eqref{eq:P1_T2r} for any $\nu \in \Psi^t$, we get
\begin{align}
    \hspace{-0.2cm} (\nu^t - \nu) \nabla_{\nu} \calL^{t} &\geq \frac{1}{2 \alpha} (|\nu^{t+1} - \nu |^2 - |\nu^t - \nu |^2) - \frac{\alpha}{2} H^2.
    \label{eq:P2_T2r}
\end{align}

The time-varying Lagrangian differences for a fixed $\bx^t \in \calX^t$ and $\blambda^t \in \Lambda^t$ satisfies \vspace{-0.1cm} 
\begin{align}
    \hspace{-0.3cm} \calL^t(\bx^t, \nu^t, \blambda^t)  - \calL^t(\bx^t, \nu, \blambda^t) &=  (\nu^t - \nu) \nabla_{\nu} \calL^t(\bx^t, \nu, \blambda). 
    \label{eq:P3_T2r}
\end{align}
Replacing \eqref{eq:P3_T2r} in \eqref{eq:P2_T2r}, and  
subtracting the results in \eqref{eq:P7_T1}, for any $\nu \in \Psi^t$, we get 
\begin{align} \nonumber
    &\calL^t(\bx^t, \nu, \blambda^t) - \calL^t(\bx^{t*}, \nu^t, \blambda^t) \\\nonumber
    &\quad \leq \frac{1}{2 \alpha} (\| \bx^t - \bx^{t*} \|^2 - \|\bx^{t+1} - \bx^{t+1*}\|^2) \\\nonumber
    &\quad \quad +  \frac{1}{2 \alpha} (|\nu^t - \nu|^2- |\nu^{t+1} - \nu |^2) + \| \be^t \| \, (2B_x + \alpha \Gamma_x) \\ 
    &\quad \quad + \frac{\alpha}{2} \| \be^t \|^2 + \frac{1}{2 \alpha} (\phi^t)^2 + \frac{1}{\alpha} D_x \phi^t + \frac{\alpha}{2} (\Gamma_x^2 + H^2).
    \label{eq:P5_T2r}
\end{align}

Since $\blambda^{t\top} \bD \bx^{t*} = 0$, we can also express $\calL^t(\bx^t, \nu, \blambda^t) - \calL^t(\bx^{t*} , \nu^t, \blambda^t)$ as
\begin{align} \nonumber
    &\calL^t(\bx^t, \nu, \blambda^t) - \calL^t(\bx^{t*} , \nu^t, \blambda^t) = \left( f(\bx^t) - f(\bx^{t*}) \right) \\ 
    & \quad \quad + \nu \,  C^t(\by^t(\bx^t)) + \blambda^{t\top} \bD \bx^t  - \nu^t C^t(\by^t(\bx^{t*})) .
    \label{eq:P1_T2a}
\end{align}
Reorganizing terms in \eqref{eq:P1_T2a}, and by the fact that $\nu^t C^t(\by^t(\bx^{t*})) \leq 0$,  we have
\begin{align} \nonumber
    \nu \, C^t(\by^t(\bx^t)) &  \leq  \left( \calL^t(\bx^t, \nu, \blambda^t) - \calL^t(\bx^{t*} , \nu^t, \blambda^t)  \right)  \\
    & - \left(f(\bx^t) - f(\bx^{t*}) \right) -  \blambda^{t \top} \bD \bx^t . 
    \label{eq:P2_T2aa}
\end{align}
On the other hand, since the cost function is convex, we have
\begin{align}
    f(\bx^t) - f(\bx^{t*}) 
    &\leq \|\bx^t - \bx^{t*}\| \| \nabla f(\bx^{t*}) \| \leq D_x L,
    \label{eq:P3_T2aaa}
\end{align}
and because $\calX^t$ and $\Lambda^t$ are compact sets uniformly in time, we further get
\begin{align}
    | \blambda^{\top} \bD \bx^t | &\leq \|\blambda^t \| \| \bD \| \|\bx^t\| \leq B_\lambda \Omega B_x.
    \label{eq:P4_T2aa}
\end{align}
By using \eqref{eq:P5_T2r}, \eqref{eq:P3_T2aaa} and \eqref{eq:P4_T2aa} in \eqref{eq:P2_T2aa}, and reorganizing terms, 
\begin{align} \nonumber
    &C^t(\by^t(\bx^t)) \leq \nu^{-1} \left( \frac{1}{\alpha} (4 B_x^2 +  B_\nu^2) + D_x L + B_\lambda \Omega B_x  \right. \\ \nonumber
    &\left. \quad \quad + \frac{\alpha}{2} (\Gamma_x^2 + H^2) + \| \be^t \| \, (2B_x + \alpha \Gamma_x) + \frac{\alpha}{2} \| \be^t \|^2 \right. \\
    &\left. \quad \quad + \frac{1}{2 \alpha} (\phi^t)^2 + \frac{1}{\alpha} D_x \phi^t \right). 
    \label{eq:P6_T2aa}
\end{align}

We now take the max operator, sum over time, and take the expectation of both sides of \eqref{eq:P6_T2aa}, then
\begin{align*} \nonumber
    &\bbE \left[ \sum_{t=1}^T \left[ C^t(\by^t(\bx^t)) \right]^+ \right] \leq T \nu^{-1} \left( \frac{1}{\alpha} (4 B_x^2 +  B_\nu^2) + D_x L  \right. \\ \nonumber
    &\left. \quad \quad + B_\lambda \Omega B_x + \frac{\alpha}{2} (\Gamma_x^2 + H^2) \right) + \nu^{-1} \left( \xi^T (2B_x + \alpha \Gamma_x) \right)  \\
    &\quad \quad + \nu^{-1} \left( \frac{\alpha}{2} \Xi^T + \frac{1}{2 \alpha} \Upsilon^T + \frac{1}{\alpha} D_x \Phi^T \right) \, \forall \, \nu \in (0, B_\nu]. 
\end{align*}
The previous inequality implies that the tightest result happens at $\nu = B_\nu$, then, \eqref{eq:P5_T2a} holds. \hfill $\blacksquare$

Finally, we present an additional result for the dynamic fit. Let $\mathrm{Fit}_T := \left[ \sum_{t=1}^T C^{t}(\by^{t}(\bx^t)) \right]^+$. By using \eqref{eq:P5_T2r}, \eqref{eq:P3_T2aaa} and \eqref{eq:P4_T2aa} in \eqref{eq:P2_T2aa}, summing over time, and using the telescopic property for the first two terms on the RHS, we get that for any $\nu \in \Psi^t$, and initial dual variable $\nu^1 = 0$

\begin{small}
\begin{align} \nonumber
    &\nu \sum_{t=1}^T C^t(\by^t(\bx^t)) \leq \frac{1}{2 \alpha} (\| \bx^1 - \bx^{1*} \|^2 + | \nu |^2) + \frac{\alpha}{2} T (\Gamma_x^2 + H^2) \\\nonumber
    &~~\quad \quad + T D_x L + T B_\lambda \Omega B_x + (2B_x + \alpha \Gamma_x) \sum_{t=1}^T \| \be^t \| \\ 
    &~~\quad \quad + \frac{\alpha}{2} \sum_{t=1}^T \| \be^t \|^2 + \frac{1}{2 \alpha} \Upsilon^T + \frac{1}{\alpha} D_x \Phi^T. 
    \label{eq:P8_T2aa}
\end{align}
\end{small}

\noindent
By taking the max operator and the expectation of both sides on \eqref{eq:P8_T2aa}, we obtain for any $\nu \in (0, B_\nu]$
\begin{align*} \nonumber
    &\bbE \left[ \left[ \sum_{t=1}^T C^t(\by^t(\bx^t)) \right]^+ \right] \leq \nu^{-1} ( T D_x L + T B_\lambda \Omega B_x ) \\\nonumber
    &~~\quad + \nu^{-1} \left( \frac{1}{2 \alpha} (4 B_x^2 + B_\nu^2) + \frac{\alpha}{2} T (\Gamma_x^2 + H^2) \right) \\ 
    &~~\quad + \nu^{-1} \left( \frac{\alpha}{2} \Xi^T + \xi^T (2B_x + \alpha \Gamma_x) + \frac{1}{2 \alpha} \Upsilon^T + \frac{1}{\alpha} D_x \Phi^T \right). 
\end{align*}
Similar as the ACV, set $\nu = B_\nu$ in order to get the tightest bound for the expected dynamic fit, i.e., $\bbE [\mathrm{Fit}_T]$.

\section{Numerical Results} \label{sec:results}

In this section, the proposed algorithm is numerically evaluated on a problem related to real-time management of distributed energy resources (DERs) at the power distribution level. In particular, we consider an aggregation of controllable DERs as well as uncontrollable loads connected to a distribution transformer. DERs include HVAC systems, electric vehicles, and energy storage systems. We consider $M=3$ DERs shared among $N=6$ users; specifically, $N_1=2$, $N_2=3$, $N_3=1$. The goal is to minimize the discomfort or dissatisfaction for each user  while enabling the aggregation of DERs to actively emulate a virtual power plant where the total power $y^t = \mathbf{1}^\top \bx_{\mathrm{in}}^t + \mathbf{1}^\top \bw^t$ follows a given automatic gain control or demand response reference signal $y_{\mathrm{ref}}^t$. The operational sets for the devices are: battery $\mathcal{X}_1 = [-8,8] \text{ kW}$; HVAC $\mathcal{X}_2 = [0, 10] \text{ kW}$; and electric vehicle $\mathcal{X}_3 = [2, 30] \text{ kW}$. We note that, in the context of inverter-interfaced DERs, the set of feasible active power setpoints is typically convex.
The discomfort functions $\{ U_{m,n} \}$ are assumed to be quadratic; the minimum of each of the functions is inside the set constraints $\mathcal{X}_m$, and it corresponds to a preferred setting of the user. For example, for EVs they represent a preferred charging rate; for HVAC systems, they represent a preferred temperature setpoint (converted into a preferred power setpoint)~\cite{Lesage_2020}. 

\begin{figure}[!ht]
  \centering 
  \begin{subfigure}[]{\includegraphics[width=0.5\textwidth]{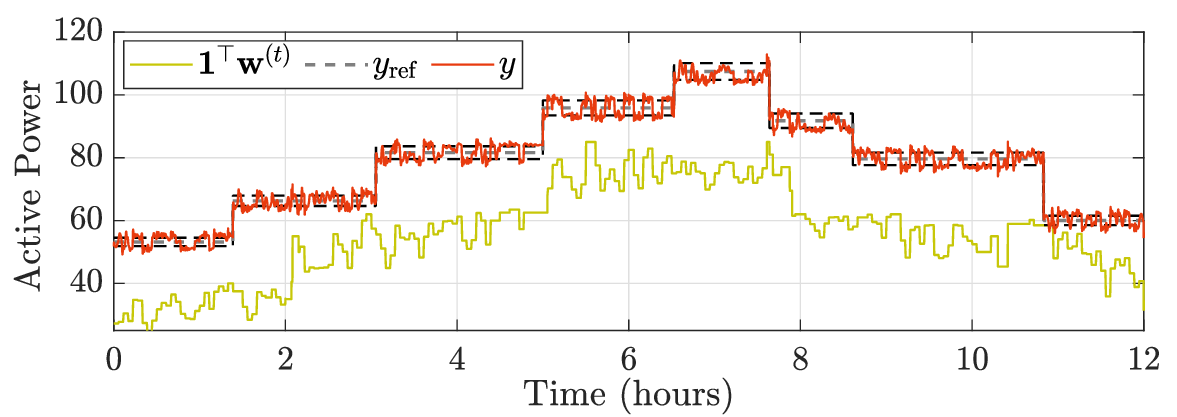}} 
  \end{subfigure}
  \begin{subfigure}[]{\includegraphics[width=0.5\textwidth]{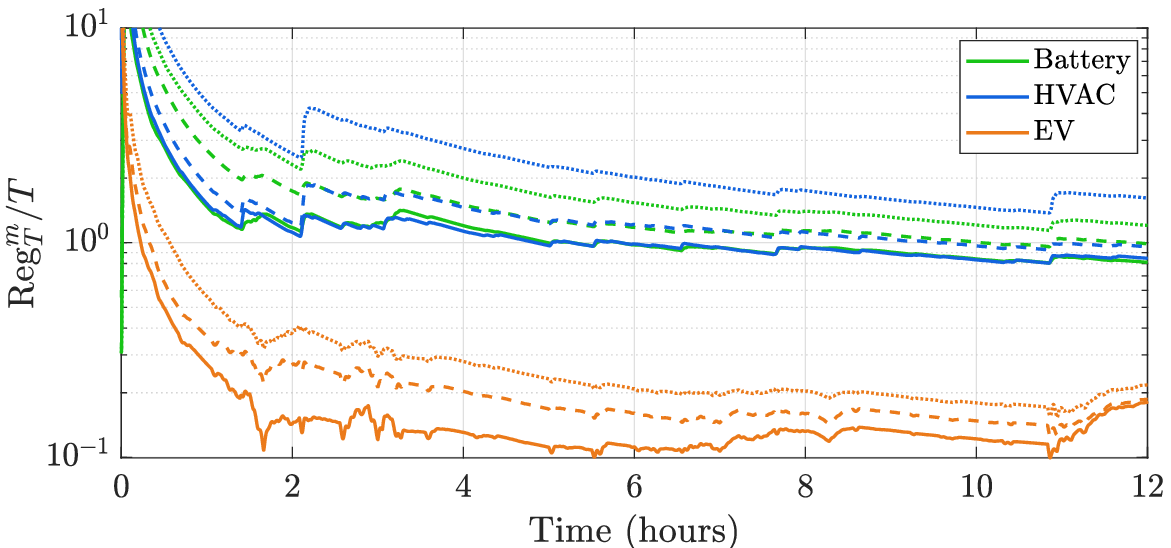}}
  \end{subfigure}
  \caption{(a) Tracking of $y_{\mathrm{ref}}$ by the GP-based online primal-dual method and (b) Local network regret (solid line corresponds to true $\{U_{m,n}\}$, dashed line to estimate $\{\hat{U}_{m,n}\}$ via shape-constrained GPs, and dotted line to estimate $\{\hat{U}_{m,n}\}$ via GPs).}
  \label{fig:regret_local}
  \vspace{-.3cm}
\end{figure}

The time-varying constraint is $C^t(\bx_\mathrm{in}^t)=\frac{\beta}{2}(\mathbf{1}^\top \bx_\mathrm{in}^t + \mathbf{1}^\top \bw^t - y_{\mathrm{ref}}^t)^2 - \zeta^t$; where $\zeta^t$ is a given tolerance, which is set to 5\% of the value of $y_{\mathrm{ref}}^t$. The power of the uncontrollable loads is  taken from the Anatolia dataset (National Renewable Energy Laboratory, Tech. Rep. NREL/TP-5500-56610), and have a granularity of 1 second (see Figure~\ref{fig:regret_local}). In this case, the evaluation of the gradient of $C^t$ requires measurements of the total power at each step $t$. 

We evaluate the performance of the online algorithm at a period of time equivalent to 12 hours; each step of Algorithm \ref{algo} is performed every 5 seconds (except for the HVAC systems, which are updated at a slower rate). The priors $\{ \hat{U}_{m,n} \}$ are determined from some noisy measurements ($\sigma = 1.5$) and the discomfort function is updated through the user's feedback every 30 min. In this case the parameter $L_U$ and $\gamma_U$ are defined beforehand; however, these values can be estimated via cross-validation. In Figure \ref{fig:regret_local}(a) it can be seen that the trajectory $y$ is within 5\% of the reference setpoint $y_\mathrm{ref}$ most of the time, in spite of the variation of the uncontrollable loads.  

\begin{figure}[!ht] 
    \centering
\includegraphics[scale=0.45]{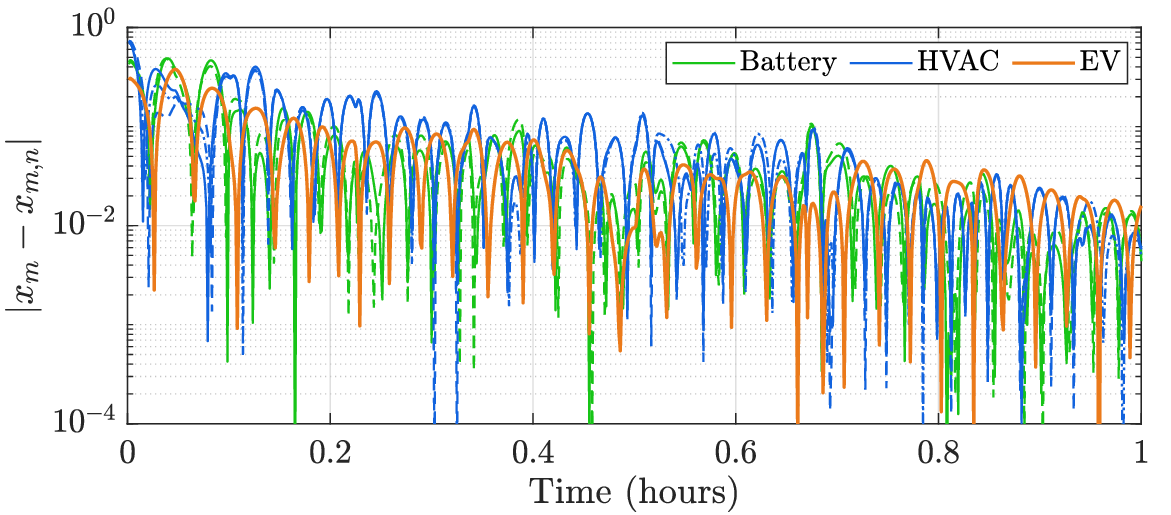}
   \caption{Disagreement between system and users (zoom for the first hour of simulation where each user is represented by a different line style). The primal variables $x_m$ and $x_{m,n}$ are normalized such that $10^{-2}$ corresponds to $1\%$ error.}
    \label{fig:disagree}
    \vspace{-.3cm}
\end{figure}

\begin{figure}[!ht] 
    \centering
\includegraphics[scale=0.45]{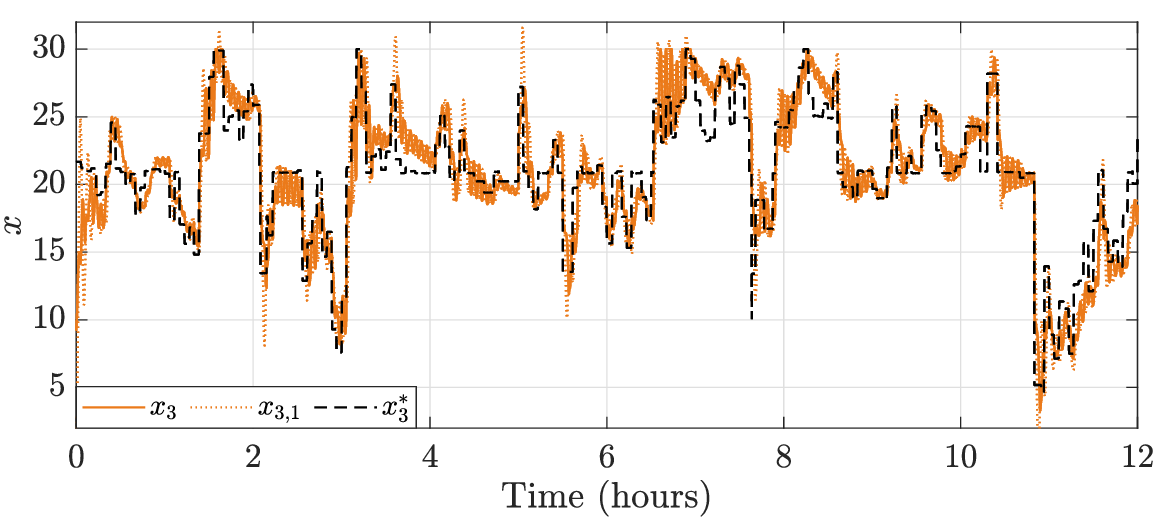}
   \caption{Behavior of the setpoint $x_3$ and the desirable setpoint for the user $x_{3,1}$, compared with the optimal solution $x_3^*$}
    \label{fig:track_devices3}
    \vspace{-.3cm}
\end{figure}

\begin{figure}[!ht] 
    \centering
\includegraphics[scale=0.45]{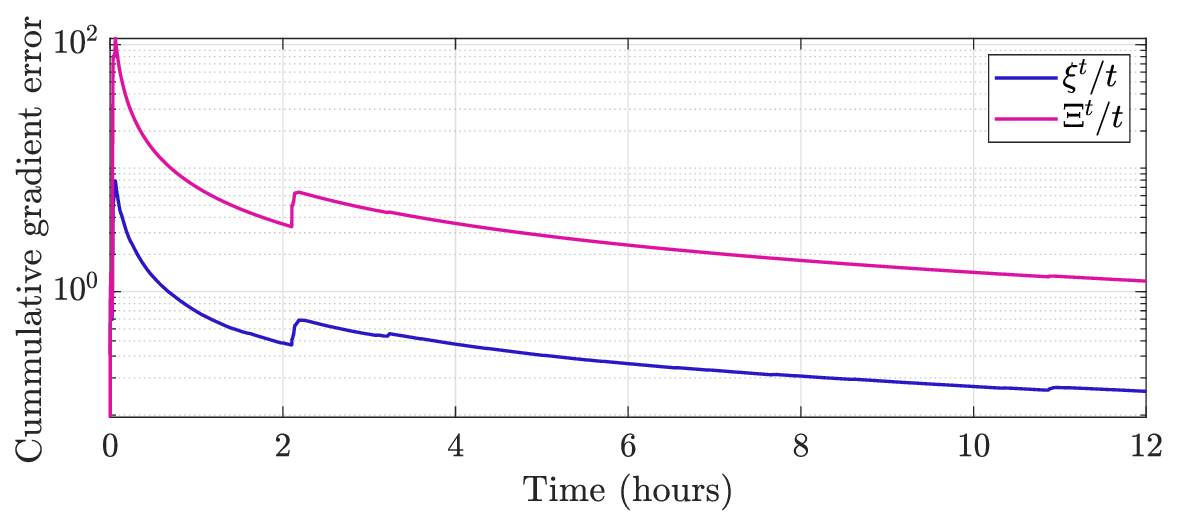}
   \caption{Behavior of $\xi^t$ and $\Xi^t$ over time.}
    \label{fig:cum_e_grad}
    \vspace{-.3cm}
\end{figure}

The local network regret is presented in Figure \ref{fig:regret_local}(b). It can be seen that the jumps in the dynamic regret corresponds to instants where the reference $y_{\text{ref}}$ changes abruptly. The regret of the proposed algorithm is, as one would expect, higher than the clairvoyant case where the functions $\{U_{m,n}\}$ are known. However, the difference diminishes over time, as the estimation error decreases. Further, Figure \ref{fig:disagree} shows the behavior of the users' disagreement during the first hour of the simulation, while Figure \ref{fig:track_devices3} illustrates the trajectory of the individual preferences of the user connected to the EV (system 3). Additionally, Figure \ref{fig:cum_e_grad} shows how the metrics in \eqref{eq:def2}, related to the accuracy of the gradient estimates of the cost function, decrease over time. 

\begin{figure}[!ht] 
    \centering
\includegraphics[scale=0.44]{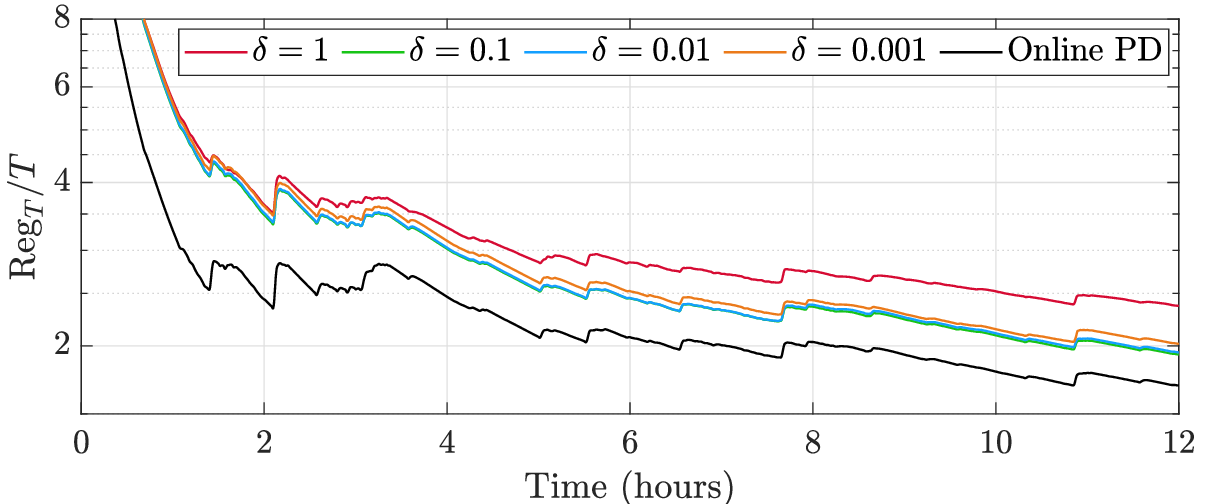}
   \caption{Global network regret for different values of $\delta$.}
    \label{fig:sen_delta}
    \vspace{-.3cm}
\end{figure}

\begin{figure}[!ht] 
    \centering
\includegraphics[scale=0.44]{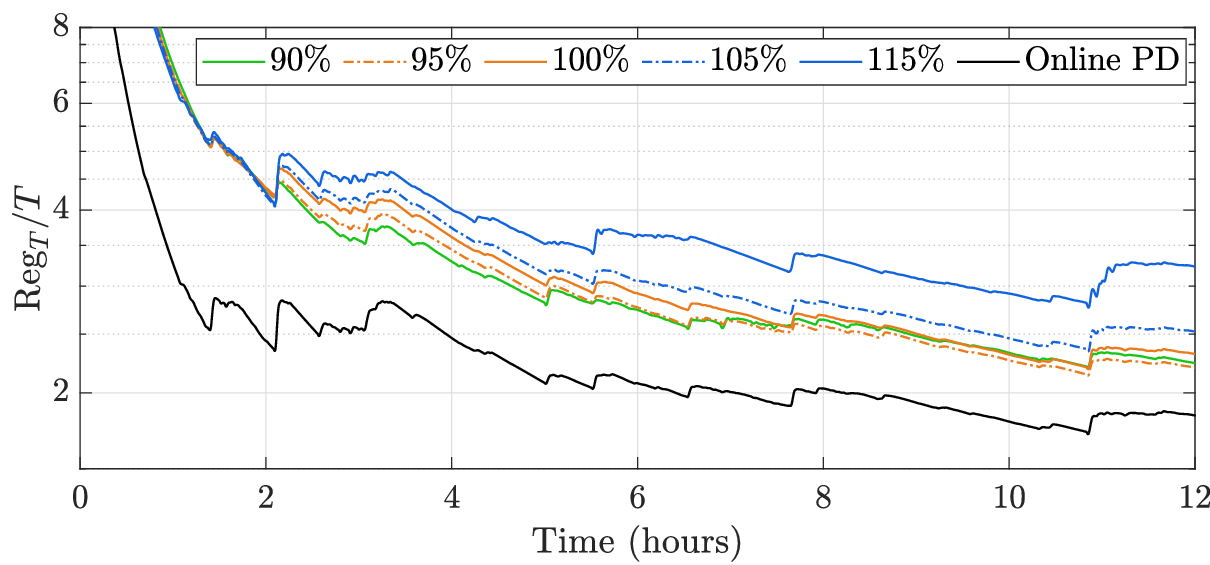}
   \caption{Global network regret over time under different load profiles.}
    \label{fig:sen_load}
    \vspace{-.3cm}
\end{figure}

\begin{figure}[!ht] 
    \centering
\includegraphics[scale=0.44]{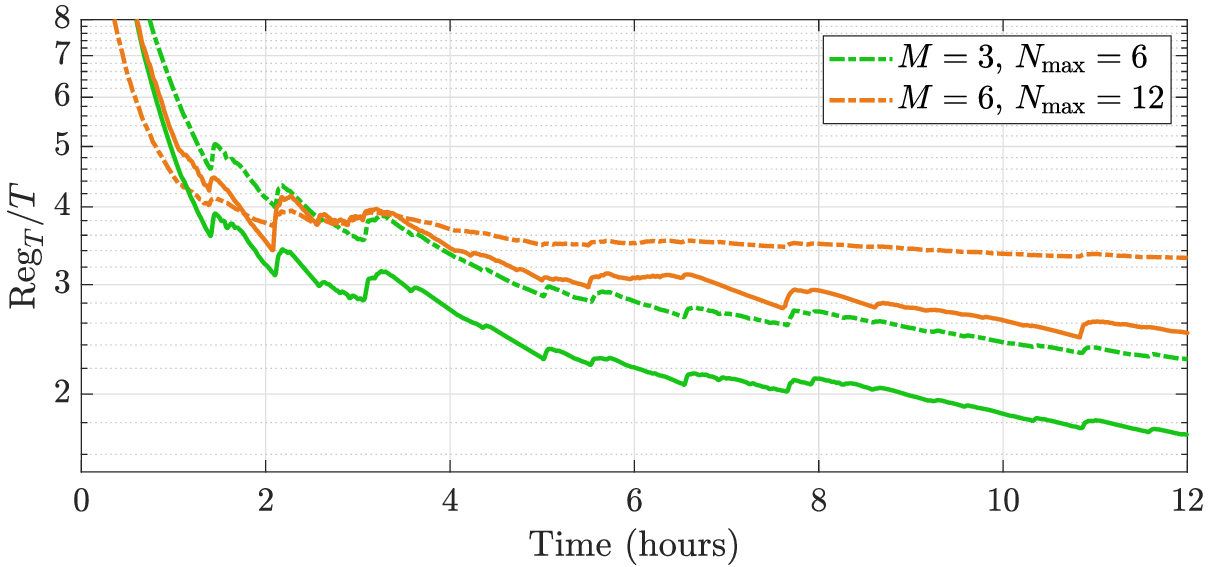}
   \caption{Global network regret for different values of $M$ and $N_\mathrm{max}$ (solid line corresponds to true $\{U_{m,n}\}$, dashed line to estimate $\{\hat{U}_{m,n}\}$ via shape-constrained GPs).}
    \label{fig:sen_M_N}
    \vspace{-.3cm}
\end{figure}

Finally, we test the performance of Algorithm \ref{algo} in terms of the global network regret under three cases. First, Figure \ref{fig:sen_delta} illustrates the sensitivity of the proposed method to the choice of the parameter $\delta$ in \eqref{eq:gradient_U}. Second, Figure \ref{fig:sen_load} presents the behavior under variability of the uncontrollable loads, in particular when the max/min load bounds increase/decrease. Lastly, Figure \ref{fig:sen_M_N} shows that the global network regret increases when the total number of devices $M$ and the total number of user of the network $N_\mathrm{max}$ increase, as expected. In these three cases, we provide the trajectory for the online algorithm with perfect knowledge of the users' function (Online PD) for comparison purposes.

\section{Conclusions}

In this work, we presented an online consensus-based algorithm to solve a time-varying optimization problem associated with a network of systems shared by multiple users. The cost function that represents the individuals' preferences were learned concurrently with the execution of the algorithm via shape-constrained GPs using users' feedback. We developed an online algorithm based on a primal-dual method, properly modified to accommodate feedback from the users and measurements from the network. We analyzed the performance via dynamic network regret and average constraint violation, and we showed that the bounds of these metrics depend on the temporal variability of the optimal solution set, and the errors associated with the estimation of the gradients. Numerical results were presented in the context of real-time management of distributed energy resources. Future efforts include investigating the convergence of multi-time-scale algorithms, and the case where the users' cost functions are time-varying to reflect preferences that change over time.


\bibliographystyle{IEEEtran}
\bibliography{References}

\end{document}